\documentclass[sn-mathphys-num]{sn-jnl}

\usepackage{amsmath}
\usepackage{hyperref}
\usepackage[
    capitalise,
    noabbrev,
  ]{cleveref}

\usepackage{amsthm}
\usepackage{graphicx}%
\usepackage{multirow}%
\usepackage{amssymb,amsfonts}%

\usepackage{mathrsfs}%
\usepackage[title]{appendix}%
\usepackage{xcolor}%
\usepackage{textcomp}%
\usepackage{manyfoot}%
\usepackage{booktabs}%
\usepackage{algorithm}%
\usepackage{algorithmicx}%
\usepackage{algpseudocode}%
\usepackage{listings}%
\usepackage{appendix} 

\usepackage{numprint} 

\npdecimalsign{.} 

\numberwithin{equation}{section}
\usepackage{xparse}

\newcommand{\R}{\mathbb R}

\newcommand{\Z}{\mathbb Z}
\newcommand{\N}{\mathbb N}

\newcommand{\Div}{\operatorname{div}}

\newcommand{\<}{\langle}
\renewcommand{\>}{\rangle}
\renewcommand{\phi}{\varphi}
\renewcommand{\epsilon}{\varepsilon}
\newcommand{\Laplace}{\operatorname{\Delta}}
\DeclareDocumentEnvironment{math}{o}{\begin{equation*}\begin{aligned}}{\end{aligned}\end{equation*}\par}

\DeclareDocumentCommand{\argdot}{}{\operatorname{\,\cdot\,}}
\DeclareDocumentCommand{\di}{O{x}}{\;\mathrm{d}#1}

\newcommand{\nnew}[1]{\textcolor{black}{#1}}
\newcommand{\bnew}[1]{\textcolor{black}{#1}}

\usepackage{comment}
\renewenvironment{comment}{\color{black}}{}

\crefname{setting}{Setting}{Settings}
\crefname{setting-s}{Setting}{Settings}
\crefname{assumption-t}{Assumption}{Assumptions}
\DeclareDocumentCommand{\eqref}{m}{\labelcref{#1}}
\crefformat{assumption-t}{Assumption~(#2#1#3)}

\crefname{appsec}{Appendix}{Appendices}

\usepackage{aliascnt}

\newtheorem{theorem}{Theorem}[section]

\newaliascnt{lemma}{theorem}
\newtheorem{lemma}[lemma]{Lemma}
\aliascntresetthe{lemma}
\crefname{lemma}{Lemma}{Lemmas}

\newaliascnt{proposition}{theorem}
\newtheorem{proposition}[proposition]{Proposition}
\aliascntresetthe{proposition}
\crefname{proposition}{Proposition}{Propositions}

\newaliascnt{remark}{theorem}

\aliascntresetthe{remark}
\crefname{remark}{Remark}{Remarks}

\newaliascnt{definition}{theorem}

\aliascntresetthe{definition}
\crefname{definition}{Definition}{Definitions}

\crefname{theorem}{Theorem}{Theorems}

\newcommand{\vecv}[1]{\vec{#1}}

\usepackage{stackrel}

\DeclareDocumentCommand{\op}{m}{\operatorname{#1}}
 \usepackage[shortlabels]{enumitem}
\usepackage{subcaption} 
\DeclareDocumentCommand{\InputData}{mO{}}{
  \def\DataPrefix{#2}
  \input{#1}}

\def\transpose{\top}
\def\Id{\operatorname{Id}}
\usepackage{todonotes}
\makeatletter
\renewcommand{\todo}[2][]{\tikzexternaldisable\@todo[#1]{#2}\tikzexternalenable}
\newcommand{\tododone}[2][]{\tikzexternaldisable\@todo[color=green!80!black,#1]{#2}\tikzexternalenable}
\makeatother

\usepackage{pgffor}

\usepackage{tikz}
\usetikzlibrary{calc}

\usepackage{pgfplotstable}
\usepackage{pgfplots}

\usepgfplotslibrary{external}
\tikzset{
  external/system call = {%
    lualatex \tikzexternalcheckshellescape%
    -halt-on-error -interaction=batchmode -jobname "\image" "\texsource"
  }
}

\definecolor{Set2-8-1}{RGB}{102,194,165}
\definecolor{Set2-8-2}{RGB}{252,141,98}
\definecolor{Set2-8-3}{RGB}{141,160,203}
\definecolor{Set2-8-4}{RGB}{231,138,195}
\definecolor{Set2-8-5}{RGB}{166,216,84}
\definecolor{Set2-8-6}{RGB}{255,217,47}
\definecolor{Set2-8-7}{RGB}{229,196,148}
\definecolor{Set2-8-8}{RGB}{179,179,179}
\pgfplotscreateplotcyclelist{mycolor}{
    {Set2-8-1},
    {Set2-8-2},
    {Set2-8-3},
    {Set2-8-4},
    {Set2-8-5},
    {Set2-8-6},
    {Set2-8-7},
    {Set2-8-8}
}

\DeclareDocumentCommand{\Code}{m}{%
  \begingroup
  \ttfamily
  \begingroup\lccode`~=`/\lowercase{\endgroup\def~}{/\discretionary{}{}{}}%
  \begingroup\lccode`~=`[\lowercase{\endgroup\def~}{[\discretionary{}{}{}}%
  \begingroup\lccode`~=`.\lowercase{\endgroup\def~}{.\discretionary{}{}{}}%
  \catcode`/=\active\catcode`[=\active\catcode`.=\active\catcode`_=\active
  \scantokens{#1\noexpand}%
  \endgroup%
}

\DeclareDocumentCommand{\PlotImageMesh}{som}{%
  \begin{tikzpicture}
    \begin{axis}[
        colormap/blackwhite,
        hide axis,
        enlargelimits=false,
        axis equal image,
        scale only axis,
        point meta min=0,
        point meta max=1,
        #2
      ]
      \addplot[patch,shader=interp]
        table[col sep=comma, point meta=\thisrow{f_xy}]
        {#3};
      \IfBooleanTF{#1}{
        \addplot[patch,shader=interp,mesh,color=red]
          table[col sep=comma, point meta=\thisrow{f_xy}]
          {#3};
      }{}
    \end{axis}
  \end{tikzpicture}%
}

\ExplSyntaxOn
\definecolor{algcomment}{HTML}{444444}
\definecolor{algcommand}{HTML}{000055}

\RequirePackage{algpseudocode}
\cs_new:Nn \mycourse_algcommand_style: { \ttfamily \color{algcommand} }

\algrenewcommand { \algorithmiccomment } [ 1 ] {
    \hfill \mbox {
        » ~ { \slshape \color { algcomment } #1 }
    }
}
\algrenewcommand { \algorithmicwhile  } { { \mycourse_algcommand_style: while  } }
\algrenewcommand { \algorithmicdo     } { { \mycourse_algcommand_style: do     } }
\algrenewcommand { \algorithmicend    } { { \mycourse_algcommand_style: end    } }
\algrenewcommand { \algorithmicfor    } { { \mycourse_algcommand_style: for    } }
\algrenewcommand { \algorithmicif     } { { \mycourse_algcommand_style: if     } }
\algrenewcommand { \algorithmicthen   } { { \mycourse_algcommand_style: then   } }
\algrenewcommand { \algorithmicelse   } { { \mycourse_algcommand_style: else   } }
\algrenewcommand { \algorithmicrepeat } { { \mycourse_algcommand_style: repeat } }
\algrenewcommand { \algorithmicuntil  } { { \mycourse_algcommand_style: until  } }
\DeclareDocumentCommand { \Return } { m } { \State { \mycourse_algcommand_style: return} ~ #1 }

\DeclareDocumentCommand { \Input }  { m } { { \small\itshape\bfseries Input: } ~ #1 }
\DeclareDocumentCommand { \Parameters }  { m } { { \small\itshape\bfseries Parameters: } ~ #1 }
\DeclareDocumentCommand { \Output } { m } { { \small\itshape\bfseries Output: } ~ #1 }
\ExplSyntaxOff

\usepackage{pdflscape}

\usepackage{csvsimple}
\usepackage{siunitx}

\begin{document}

\title[A primal-dual adaptive finite element method for total variation minimization]{A primal-dual adaptive finite element method for total variation minimization}


\author[1]{\fnm{Martin} \sur{Alk\"amper}}

\author[2]{\fnm{Stephan} \sur{Hilb}}

\author*[3]{\fnm{Andreas} \sur{Langer}}\email{andreas.langer@math.lth.se}

\affil[1]{\city{Ratingen}, \country{Germany}}
\affil[2]{\city{Stuttgart}, \country{Germany}}
\affil*[3]{\orgdiv{Department for Mathematical Sciences}, \orgname{Lund University}, \orgaddress{\street{S\"olvegatan 18A}, \city{Lund}, \postcode{22100}, \country{Sweden}}}



\abstract{Based on previous work we extend a primal-dual semi-smooth Newton method for minimizing a general $L^1$-$L^2$-$TV$ functional over the space of functions of bounded variations by adaptivity in a finite element setting. For automatically generating an adaptive grid we introduce indicators based on a-posteriori error estimates. Further we discuss data interpolation methods on unstructured grids in the context of image processing and present a pixel-based interpolation method. 
The efficiency of our derived adaptive finite element scheme is demonstrated on image inpainting and the task of computing the optical flow in image sequences. In particular, for optical flow estimation we derive an adaptive finite element
coarse-to-fine scheme which allows resolving large displacements and speeds-up the computing time significantly.}




\maketitle

\section{Introduction}\label{sec1}

Let $\Omega \subset \R^d$, $d \in \N$, be an open, bounded and simply connected domain with Lipschitz boundary.
We consider the problem of minimizing the 
$L^1$-$L^2$-TV functional \cite{HintermullerLanger:13}, which has various applications in image processing and has been subject of study in \cite{Langer:17a,Langer:19} and in a finite element context in \cite{AlkamperLanger:17}. More precisely, \nnew{let $V \subset H^1(\Omega)^m$, $m\in\N$, be a continuously embedded Hilbert space which is weakly closed in $H^1(\Omega)^m$ and $L^2(\Omega)^m$, then as in \cite{ours2022semismooth1} we consider the $L^1$-$L^2$-TV functional} in the following form 
\begin{equation} \label{eq:primal_smoothed}
  \begin{aligned}
    \inf_{\vecv{u} \in V} \Big\{
    \alpha_1 \int_\Omega \phi_{\gamma_1}(|T \vecv{u} - g|) \di[\vecv{x}]
    &+ \tfrac{\alpha_2}{2} \|T \vecv{u} -g \|_{L^2}^2
    + \tfrac{\beta}{2} \|S\vecv{u}\|_{V_S}^2\\
    &+ \lambda \int_\Omega \phi_{\gamma_2}(|\nnew{\nabla} \vecv{u}|_F) \di[\vecv{x}]
    =: E(\vecv{u})
    \Big\},
  \end{aligned}
\end{equation}
where 
$S: V \to V_S$ is a bounded linear operator for some Hilbert space $V_S$,
$T: L^2(\Omega)^m \to L^2(\Omega)$ is a bounded linear operator,
$g \in L^2(\Omega)$ is some given data, $|\argdot|_F:\R^{d\times m} \mapsto \R$ denotes the Frobenius norm, and
$\alpha_1, \alpha_2, \lambda, \beta \ge 0$ are adjustable weighting parameters. Further
$\gamma_1, \gamma_2 \ge 0$ are Huber-regularization parameters associated to
the \emph{Huber-function} $\phi_\gamma: \R \to [0,\infty)$ for $\gamma \ge 0$ defined by
\begin{equation} \label{eq:phidef}
  \phi_\gamma(x) := \begin{cases}
    \frac{1}{2\gamma} x^2 & \text{if } |x| \le \gamma, \\
    |x| - \frac{\gamma}{2} & \text{if } |x| > \gamma.
  \end{cases}
\end{equation}
We remark that $d$ denotes the spatial dimension, e.g.\ one sets $d=1$ for signals and $d=2$ for images, and $m$ describes the number of channels, i.e.\ for grey-scale images $m=1$ and for motion fields we have $m=d$. 
In \cite{ours2022semismooth1} a primal-dual framework around \eqref{eq:primal_smoothed} has been developed and a finite element semi-smooth Newton scheme established. The present work aims to extend \cite{ours2022semismooth1} to an adaptive finite element setting. 
\nnew{Note that every finite dimensional subspace of a normed vector space is closed \cite[Corollary 5.34]{HunterNachtergaele:01} and convex and hence also weakly closed \cite[Corollary 8.74]{EinsiedlerWard:17}. Thus $V$ in \eqref{eq:primal_smoothed} includes any arbitrary finite element space.}
Moreover, we note that the finite element semi-smooth Newton scheme in \cite{ours2022semismooth1} is also derived for this space setting. 

\subsubsection*{Related work}

While mathematical image processing techniques are often presented or discretized in a finite difference setting due to simple handling of image array data, finite element methods in image processing have also been studied, see \cite{ChambollePock:21} for an overview. 
The ability to adaptively refine or coarsen unstructured finite element meshes
seems promising in order to save computational resources during processing. We refer the reader to \cite{NoSiVe:09,Verfurth:13} for an overview of adaptive finite element methods. 
We note that images usually have homogeneous regions as well as parts with a lot of details. It seems intuitively sufficient to discretize finely in regions of transitions and fine details, while a coarse discretization should do well in uniform parts. 
To this end error estimates are crucial in order to correctly guide the refinement process. First attempts for using adaptive finite element methods in image processing (restoration) have been made for non-linear diffusion models, which are mainly used for image denoising and edge detection. In particular, in \cite{BanschMikula:97,BazanBlomgren:05,PreuerRumpf:00} adaptive finite element methods are proposed for the modified Perona-Malik model \cite{CaLiMoCo:92}. 
Based on primal-dual gap error estimators adaptive finite element methods for the $L^2$-TV model, i.e.\ \eqref{eq:primal_smoothed} with $\alpha_1=\beta=\gamma_2=0$, are proposed in \cite{Bartels:15a,BartelsKaltenbach:24,BartelsMilicevic:20}. 
While in \cite{Bartels:15a} the dual problem is discretized by continuous, piecewise affine finite elements, in \cite{BartelsMilicevic:20} the Brezzi-Douglas-Morini finite elements \cite{BrDoMa:85} are used. 
\nnew{Additionally, \cite{BartelsKaltenbach:24} introduces an adaptive finite element framework based on Crouzeix–Raviart elements.}
Note that in \cite{Bartels:15a,BartelsKaltenbach:24,BartelsMilicevic:20} only clean images are considered and no image reconstruction or analysis problem is tackled. In a similar way relying on a semi-smooth Newton method for minimizing the $L^2$-TV model over the Hilbert space $H^1_0(\Omega)$ in \cite{HintermullerRincon-Camacho:14} an adaptive finite element method, again using a posteriori error estimates, is suggested for image denoising only. 
Adaptive finite element methods for optical flow computation have been presented in \cite{belhachmi2011control,BelHec2016}. However, these methods differ significantly from our approaches. Firstly, a different minimization problem is considered. Secondly, the mesh in the methods in \cite{belhachmi2011control,BelHec2016} is initialized at fine image resolution and iteratively coarsened only after a costly computation, while our approach allows to start with a coarse mesh and refines cells only if deemed necessary.

\subsubsection*{On this work}

In the present work we extend the finite element framework of \cite{ours2022semismooth1} to an adaptive finite element setting and demonstrate its efficiency by considering the applications of inpainting and motion estimation. 
The proposed adaptive finite element method is based on a-posteriori error estimates. %
\begin{comment}
In particular, we derive two different a-posteriori error estimators, one using the residual and the other using a primal-dual approach.
In order to obtain a primal-dual gap error estimate we use the same approach as in \cite{BartelsMilicevic:20}. However, due to the presence of Huber regularization terms and potentially non-local operators $T$ in \eqref{eq:primal_smoothed} we need to adjust arguments of \cite{BartelsMilicevic:20}.
\end{comment}
For deriving the residual-based error estimates we do not consider a variation inequality setting as in \cite{HintermullerRincon-Camacho:14} but utilize the strong convexity of the functional in \eqref{eq:primal_smoothed}. Further, while the approach in \cite{HintermullerRincon-Camacho:14} only considers the image denoising case, i.e.\  $T$ being the identity, our estimate allows for general linear and bounded operators $T$. 
\nnew{Utilizing our derived a-posteriori estimators allows us to propose an adaptive finite element scheme for solving \eqref{eq:primal_smoothed} for general bounded and linear operators $T$.} 
\bnew{Note that implementing global operators on an unstructured grid is generally considered impractical due to the inefficiency caused by the non-local access pattern. Therefore, we restrict our numerical experiments in this paper to local operators $T$.}

We demonstrate the usability of the adaptive scheme by solving the task of image inpainting and for estimating motion in image sequences. For both applications an adaptive coarse-to-fine scheme which makes use of the developed adaptive finite element framework is used. That is, the algorithm starts with a coarse mesh, typically coarser than the observed image, and adaptively refines the mesh at parts where deemed necessary. 
In the setting of optical flow computation this means that we incorporate into \cite[Algorithm 3]{ours2022semismooth1} an adaptive coarse-to-fine warping scheme which is able to resolve large displacements. In particular, due to the adaptive coarse-to-fine warping scheme the proposed algorithm improves \cite[Algorithm 3]{ours2022semismooth1} in two directions. Firstly, it produces more accurate approximations of the optical flow and secondly, it significantly reduces the computing time, see \cref{sec:numerical-examples} below.

Note that for finite elements in general an interpolation method needs to be employed to
transfer given data onto a finite element function. The choice of this method is especially relevant when the underlying grid is coarser than the original image resolution since in that case e.g.\ nodal sampling for standard linear Lagrange finite elements can lead to aliasing artifacts. 
In this vein, we propose a new pixel-adapted method to interpolate image data onto
finite element meshes and compare it to other alternatives. This new method is designed to minimize the discrete $\ell^2$-distance to the original image on the original image grid. That is, it actually  optimizes the peak signal-to-noise ratio ($\op{PSNR}$) and therefore might be especially suited for images.

The rest of the paper is organized as follows: In \cref{sec:Preliminaries} we fix the basic notation and terminology used and describe the primal-dual setting of \eqref{eq:primal_smoothed}. Based on this primal-dual formulation we introduce in \cref{sec:finite-elements} our finite element setting. The main part of this section is devoted to the derivative of a-posteriori error estimates to establish an adaptive finite element algorithm for \eqref{eq:primal_smoothed}. Moreover, we describe several different image interpolation methods and compare them with a newly proposed pixel-adapted interpolation method. Based on the developed adaptive finite element method in \cref{sec:numerical-examples} adaptive coarse-to-fine schemes for image inpainting and motion estimation are proposed. Numerical experiments demonstrate the applicability of the proposed adaptive discretization method. 

\section{Preliminaries}\label{sec:Preliminaries}
\subsection{Notation and Terminology}
For a Hilbert space $H$ we denote by $\<\argdot, \argdot\>_H$ its inner product and by $\|\argdot\|_H$ its induced norm. The duality pairing between $H$ and its continuous dual space $H^*$ is written as $\<\argdot, \argdot\>_{H, H^*}$. In the case $H=L^2(\Omega)^{n}$, $n \in \N$, we define the associated inner product by $\<\argdot, \argdot\>_H: \big((u_k)_{k=1}^n, (v_k)_{k=1}^n\big) \mapsto \sum_{k=1}^n \<u_{k}, v_{k}\>_{L^2(\Omega)}$, where $\<\argdot, \argdot\>_{L^2(\Omega)}$ is the standard $L^2$ inner product. Apart from notational convenience, we treat a matrix-valued space $L^2(\Omega)^{d\times m}$ as equivalent to $L^2(\Omega)^{dm}$ using
the numbering $(i,j) \mapsto (j-1) d + i$, $i \in \{1, \dotsc, d\}$, $j \in \{1, \dotsc, m\}$ of the respective components. Moreover, for any $L^2$ function space we may use the inner product shorthand notations $\<\argdot, \argdot\>_{L^2} := \<\argdot, \argdot\>_H$ and similarly $\|\argdot\|_{L^2} := \|\argdot\|_H$ for the norm.
In the sequel, the dual space $H^*$ of a Hilbert space $H$ will often be identified with itself. 

For a bounded linear operator $\Lambda\in \mathcal{L}(X,Y)$, where $X,Y$ are Hilbert spaces, we define by $\|\Lambda\| := \|\Lambda\|_{\mathcal{L}(X,Y)}$ its operator norm and denote by $\Lambda^*:Y^* \to X^*$ its adjoint operator. In the sequel by single strokes, i.e.\ $|\argdot|$, we describe the Euclidean norm on $\R^n$, $n\in\N$. Often operations are applied in a pointwise sense, such that for a vector-valued function $\vecv{u}: \Omega \to \R^m$ the expression $|\vecv{u}|$ denotes the
function $|\vecv{u}|: \Omega \to \R$, $\vecv{x} \mapsto |\vecv{u}(\vecv{x})|$.

A function $f: H \to \overline{\R}:=\R \cup\{+\infty\}$ is called \emph{coercive}, if for any sequence $(\vecv{v}_n)_{n\in\N} \subset H$ we have
  \begin{align*}
    \|\vecv{v}_n\|_H \to \infty
    \implies
    f(\vecv{v}_n) \to \infty.
  \end{align*}
A bilinear form $a:H\times H \to\R$ is called \emph{coercive}, if there exists a constant $c>0$ such that $a(\vecv{v},\vecv{v}) \ge c \|\vecv{v}\|^2_H$ for all $\vecv{v}\in H$.
The subdifferential for a convex function $f: X \to \overline{\R}:=\R \cup \{+\infty\}$ at $\vecv{v}\in X$ is defined as the set valued function
\begin{equation*}
\partial f(\vecv{v}) = \begin{cases}
\emptyset & \text{if } f(\vecv{v}) = \infty,\\
\{\vecv{v}^* \in X^*  :  \langle \vecv{v}^*, \vecv{u}-\vecv{v}\rangle_{X^*,X} + f(\vecv{v}) \leq f(\vecv{u}) \  \forall \vecv{u}\in X \} & \text{otherwise}.
\end{cases}
\end{equation*} 
Further for a predicate $w: \Omega \to \{\mathrm{true}, \mathrm{false}\}$ we define the indicator $\chi_w \in \overline{\R}$ as
\begin{equation*}
  \chi_{w} := \begin{cases}
    0 & \text{if $w(\vecv{x})$ is true for almost every\ $\vecv{x} \in \Omega$}, \\
    \infty & \text{else}.
  \end{cases}
\end{equation*}
For example, let $u \in L^2(\Omega)$, then $\chi_{|u| \le 1}$ would evaluate to $\infty$ if and only if $|u|$ is greater than $1$ on a set of non-zero measure.

\subsection{Primal-Dual Formulation}
\label{ssec:predual-and-dualization}
Let us recall the functional setting and the primal-dual formulation of \cite{ours2022semismooth1}, on which our derived a-posteriori estimates and hence our proposed adaptive finite element scheme, presented in \cref{sec:finite-elements} below, are based. 
We restrict ourselves for the operator $S$  in \eqref{eq:primal_smoothed} and its related spaces to the following two settings:
\begin{enumerate}[leftmargin=*,align=left,label=($S$.\roman*)]
  \item \label[setting-s]{setting-identity}
    $S := I: V \to V_S$ with 
    $\|\argdot\|_V := \|\argdot\|_{\nnew{L^2}}$ 
    and
    $V_S = L^2(\Omega)^m$, $\|\argdot\|_{V_S}:= \|\argdot\|_{L^2}$ or 
  \item \label[setting-s]{setting-nabla}
    $S := \nabla: V \to V_S$ with
    $\|\argdot\|_V := \|\argdot\|_{H^1}$ and
    $V_S = L^2(\Omega)^{d\times m}$ 
    (the boundedness of $S$ follows due
    to $\|\nabla \vecv{v}\|_{L^2} \leq \|\vecv{v}\|_{H^1}$),
    $\|\argdot\|_{V_S}:= \|\argdot\|_{L^2}$.
\end{enumerate}

We define the symmetric bilinear form $a_B :V\times V \to \R$ by
\begin{equation} \label{eq:bilin}
  a_B(\vecv{u},\vecv{w})
  := \alpha_2 \<T \vecv{u}, T \vecv{w}\>_{L^2} + \beta \<S\vecv{u}, S\vecv{w}\>_{L^2}
  = \<B\vecv{u}, \vecv{w}\>_{V^*,V},
\end{equation}
whereby $B := \alpha_2 T^* T + \beta S^* S: V \to V^*$. The bilinear form $a_B$ induces a respective energy norm defined
by $$\|\vecv{u} \|_B^2 := a_B(\vecv{u},\vecv{u}) \qquad \text{for} \ \vecv{u}\in V.$$

To guarantee the existence of a (unique) solution of \eqref{eq:primal_smoothed} in the sequel we assume that $a_B$ is coercive \cite[Proposition 3.1]{ours2022semismooth1}:
\makeatletter
\newcommand{\leqnomode}{\tagsleft@true\let\veqno\@@leqno}
\makeatother
\begin{equation}\leqnomode
\tag{As1}\text{The bilinear form $a_B: V \times V \to \R$ is coercive.} \label[assumption-t]{assumption-a-coercivity}
\end{equation}
Here the coercivity is to be understood in the following ways according to \cref{setting-identity} and \cref{setting-nabla}: In \cref{setting-identity} the coercivity is with respect to the $L^2$-norm, i.e.\ $a_B$ is $L^2$-coercive. \nnew{For \cref{setting-nabla} the coercivity is understood with respect to $\|\cdot\|_V$ (i.e.\ the $H^1$-norm). Since $V$ is weakly closed in $H^1(\Omega)^m$ as well as $L^2(\Omega)^m$, the existence of a solution in $V$ is guaranteed in both settings \cite[Proposition 2.2]{ours2022semismooth1}. If $V$ is finite dimensional, then \cref{assumption-a-coercivity} guarantees a solution in the respective finite dimensional space.}

Further, the coercivity of $a_B$ implies the invertibility of $B$ \cite[Remark 2.1]{ours2022semismooth1}. Accordingly, we define the norm $\|\vecv{u}\|_{B^{-1}}^2 := \<\vecv{u},
{B^{-1}}\vecv{u}\>_{V^*,V}$ for $\vecv{u}\in V^*$. Additionally, we introduce the linear and bounded operator $\Lambda := (T, \nabla): V \to W := L^2(\Omega)\times L^2(\Omega)^{d\times m}$.

In these settings we recall the following duality result.
\begin{proposition}[{\cite[Theorem 3.1 $\&$ Proposition 3.2]{ours2022semismooth1}}] \label{thm:optcond_reg_nabla}
  The pre-dual problem to \eqref{eq:primal_smoothed} reads
  \begin{equation} \label{reg_nabla_dualP}
    \begin{aligned}
      \inf_{\vec{p} = (p_1, \vec{p}_2) \in W^*}\;
      \Big\{\tfrac{1}{2} \big\|\Lambda^* \vec{p} - &\alpha_2 T^*g
        \big\|_{B^{-1}}^2 - \tfrac{\alpha_2}{2} \|g\|_{L^2}^2
        + \<g, p_1\>_{L^2}
        + \chi_{|p_1| \le \alpha_1}\\
        &+ \tfrac{\gamma_1}{2\alpha_1} \|p_1\|_{L^2}^2 
        + \chi_{|\vec{p}_2|_F\le \lambda}
        + \tfrac{\gamma_2}{2\lambda} \|\vec{p}_2\|_{L^2}^2
      =:D(\vec{p}) \Big\},
    \end{aligned}
  \end{equation}
  where for $\alpha_1 = 0$ or $\lambda = 0$ we use the convention that the terms
  $\tfrac{\gamma_1}{2\alpha_1} \|p_1\|_{L^2}^2$ and $\tfrac{\gamma_2}{2\lambda}
  \|\vec{p}_2\|_{L^2}^2$ vanish respectively.

  Furthermore, solutions $\vec{p} = (p_1, \vec{p}_2) \in W^*$ and $\vecv{u} \in V$ of
  \eqref{reg_nabla_dualP} and \eqref{eq:primal_smoothed} respectively are
  characterized by
  \begin{equation} \label{eq:optimality_condition_h1}
    \begin{aligned}
      0 &= \Lambda^* \vec{p} - \alpha_2 T^* g + B\vecv{u}, \\
      0 &= p_1 \max\{\gamma_1,|T \vecv{u} -g|\} - \alpha_1 (T \vecv{u} -g), &
        |p_1| &\le \alpha_1, \\
      0 &= \vec{p}_2 \max\{\gamma_2,|\nabla \vecv{u}|_{F}\} - \lambda \nabla \vecv{u}, &
        |\vec{p}_2|_F &\le \lambda,
    \end{aligned}
  \end{equation}
  where $\max$ denotes the pointwise maximum.
\end{proposition}

The primal functional $E$ from \eqref{eq:primal_smoothed} satisfies the
following strong convexity property.

\begin{lemma} \label{lem:energy-strong-convexity}
  If ${\vecv{u}} \in V$ is a minimizer of $E$, then we have
\begin{align*}
  \tfrac{1}{2}\|\vecv{v} - \vecv{u}\|_B^2 \leq E(\vecv{v}) - E({\vecv{u}})
\end{align*}
for all $\vecv{v}\in V.$
  \begin{proof}
    We apply the same method as in \cite[Lemma 10.2]{Bartels:12}.
    For the energy $E$ from \eqref{eq:primal_smoothed} we write
    $E(\vecv{v}) = F(\vecv{v}) + G(\vecv{v})$ with
    \begin{align*}
      G(\vecv{v}):=
        \alpha_1 \int_\Omega \phi_{\gamma_1}(|T \vecv{v} - g|) \di[\vecv{x}]
        + \lambda \int_\Omega \phi_{\gamma_2}(|\nnew{\nabla} \vecv{v}|_F) \di[\vecv{x}]
    \end{align*}
    and 
    \begin{align*}
      F(\vecv{v}):=\tfrac{\alpha_2}{2} \|T \vecv{v} -g \|_{L^2}^2 + \tfrac{\beta}{2} \|S\vecv{v}\|_{L^2}^2.
    \end{align*}
    Note that $F$ is Frechet-differentiable with
    \[
      \<F'(\vecv{v}), \vecv{w}\>_{V^*,V} = \alpha_2 \langle T \vecv{v} - g, T \vecv{w} \rangle_{L^2} + \beta \langle S\vecv{v}, S\vecv{w} \rangle_{L^2}
    \]
    for all $\vecv{w}\in V$.
    Expanding $F(\vecv{v})$ at $\vecv{u}$ then yields
    \[
      F(\vecv{v}) =
      F(\vecv{u}) + \<F'(\vecv{u}), \vecv{v} - \vecv{u}\>_{V^*,V} + \tfrac{1}{2}\|\vecv{v} - \vecv{u}\|_B^2.
    \]
    Since $\vecv{u}\in V$ is a minimizer, we have
    $0\in\partial E(\vecv{u}) = F'(\vecv{u}) + \partial G(\vecv{u})$.
    Hence $- F'(\vecv{u}) \in \partial G(\vecv{u})$, i.e.\
    \begin{align*}
      \<-F'(\vecv{u}), \vecv{v}-\vecv{u}\>_{V^*,V}
        \leq G(\vecv{v}) - G(\vecv{u}) 
    \end{align*}
    and thus 
    \begin{align*}
        \tfrac{1}{2}\|\vecv{v} - \vecv{u}\|_B^2
        \le F(\vecv{v}) -  F(\vecv{u}) + G(\vecv{v}) - G(\vecv{u}),
    \end{align*}
    which proves the statement.
  \end{proof}
\end{lemma}

Putting \cref{lem:energy-strong-convexity} and coercivity of $a_B$ from
\cref{assumption-a-coercivity}
together we obtain that for a minimizer $\vecv{u}$ of $E$ there is a constant $c > 0$
such that
\[
  E(\vecv{v}) - E(\vecv{u})
  \geq \tfrac{1}{2}\|\vecv{u} - \vecv{v}\|_B^2
  =\begin{cases}
  \tfrac{c}{2} \|\vecv{u}-\vecv{v}\|_{L^2}^2 &\text{for \cref{setting-identity}} \\
  \tfrac{c}{2} \|\vecv{u}-\vecv{v}\|_{H^1}^2 &\text{for \cref{setting-nabla}}
  \end{cases}
\]
for all $\vecv{v}\in V$.

\section{Discretization: Finite Elements}
\label{sec:finite-elements}

Central to the idea of finite element discretization is the idea to find a
solution within a finite dimensional subspace of the original space.
This discrete space usually consists of cellwise polynomial functions defined
on a mesh of cells, e.g.\ simplices.
In contrast to finite differences the discretization may be adaptively refined
to accomodate certain local error indicators, thereby reducing the amount of
degrees of freedom needed to represent a solution up to a certain accuracy. 
For more details on adaptive finite element schemes in general, we refer the
reader to the survey \cite{NoSiVe:09}.

We are interested in discretizing two-dimensional computer images of size $n_1\times n_2$ pixels, where $n_1, n_2 \in \N$, via finite elements. That is here and in the sequel we have $d=2$. For such an image the domain $\Omega$ is defined as $\Omega := [1, n_1] \times [1, n_2]$. We triangulate $\Omega$ by simplices with nodes at
integer coordinates $(x_1, x_2) \in \Z^2$, $1 \le x_1 \le n_1$, $1 \le x_2 \le
n_2$ corresponding to pixel centers as depicted in
\cref{fig:grid-discretization}.

\begin{figure}[ht]
  \centering
  \def\gridsize{0.3}
  \begin{subfigure}[t]{\gridsize\columnwidth}
    \centering
    \begin{tikzpicture}
      \foreach \x in {0, 1, 2} {
        \foreach \y [evaluate = \x as \shade using (1+\x)*(1+\y)/9*50]
            in {0, 1, 2} {
          \fill[black!\shade] (\x, \y) rectangle ++(1, 1);
        }
      }
    \end{tikzpicture}
    \caption{image with $3\times 3$ square pixels}
    \label{fig:fem-pixel-image}
  \end{subfigure}%
  \hspace{2em}%
  \begin{subfigure}[t]{\gridsize\columnwidth}
    \centering
    \begin{tikzpicture}
      \foreach \x in {0, 1, 2} {
        \foreach \y [evaluate = \x as \shade using (1+\x)*(1+\y)/9*50]
            in {0, 1, 2} {
          \fill[black!\shade] (\x, \y) rectangle ++(1, 1);
        }
      }
      \foreach \x in {0.5, 1.5} {
        \foreach \y in {0.5, 1.5} {
          \draw[black] (\x, \y) rectangle ++(1, 1);
        }
      }
      \draw (0.5, 0.5) -- (1.5, 1.5);
      \draw (2.5, 2.5) -- (1.5, 1.5);
      \draw (0.5, 2.5) -- (1.5, 1.5);
      \draw (2.5, 0.5) -- (1.5, 1.5);
      \foreach \x in {0, 1, 2} {
        \foreach \y in {0, 1, 2} {
          \filldraw ($(\x, \y) + (0.5, 0.5)$) circle [radius = 1.5pt];
        }
      }
    \end{tikzpicture}
    \caption{simplicial grid whose nodes are located in the pixel centers}
    \label{fig:fem-dual-grid}
  \end{subfigure}%
  \caption{Image aligned simplicial grid construction as in \cite{ours2022semismooth1}.}
  \label{fig:grid-discretization}
\end{figure}

Let $\mathcal T$ denote the set of cells and $\Gamma$ the set of
oriented facets (i.e.\ edges as $d = 2$) of the simplicial triangulation.
For any cell $K \in \mathcal T$ let $P_k(K)$ be the space of polynomial
functions on $K$ with total degree $k \in \N$.
We define the finite subspaces $V_h \nnew{ \subseteq V \subseteq } H^1(\Omega)^m$, $W_h^*
\subset L^2(\Omega) \times L^2(\Omega)^{d \times m} \simeq W^*$, i.e.\ $W^*$ can be identified with  $W=L^2(\Omega) \times L^2(\Omega)^{d \times m}$, and $Z_h \subset
L^2(\Omega)$ as 
\begin{equation} \label{eq:fem-spaces}
  \begin{aligned}
    V_h &:= \{\vecv{u} \in C(\Omega)^m : \vecv{u}|_K \in P_1(K)^m, K \in \mathcal T\}, \\
    W_h^* &:=
      \{(p_1, \vec{p}_2) \in C(\Omega) \times L^2(\Omega)^{d \times m}:
        p_1|_K \in P_1(K),
        \vec{p}_2|_{K} \in P_0(K)^{d \times m},\\
        &\phantom{:=
      \{(p_1, \vec{p}_2) \in C(\Omega) \times L^2(\Omega)^{d \times m}:
        p_1|_K \in P_1(K),
        \vec{p}_2|_{K} \in P_0(K)} K \in \mathcal T \}, \\
    Z_h &:= \{g \in C(\Omega) : g|_K \in P_1(K), K \in \mathcal T \}.
  \end{aligned}
\end{equation}
As in \cite{ours2022semismooth1} we do not discretize \eqref{eq:primal_smoothed} directly but the optimality conditions \eqref{eq:optimality_condition_h1} instead. In our discrete finite element setting this means that we search for
solutions $\vec{p}_{h} = (p_{\nnew{h,1}}, \vec{p}_{\nnew{h,2}}) \in W_h^*$, $\vecv{u}_h \in V_h$ which \nnew{pointwise} satisfy
\begin{equation} \label{eq:optimality_condition_h1_discrete}
  \begin{aligned}
    0 &= \Lambda^* \vec{p}_{h} - \alpha_2 T^* g_h + B\vecv{u}_{h}, \\
    0 &= p_{\nnew{h,1}} \max\{\gamma_1,|T \vecv{u}_{h} -g_h|\} - \alpha_1 (T \vecv{u}_{h} -g_h), &
      |p_{\nnew{h,1}}| &\le \alpha_1, \\
    0 &= \vec{p}_{\nnew{h,2}} \max\{\gamma_2,|\nabla \vecv{u}_{h}|_{F}\} - \lambda \nabla \vecv{u}_{h}, &
      |\vec{p}_{\nnew{h,2}}|_F &\le \lambda,
  \end{aligned}
\end{equation}
where $g_h\in Z_h$ is a discrete version of the data $g$. We elaborate further in \cref{sec:ImageInterpolation} below on possibilities how to map image data onto the given grid. Note that the last two equations in \eqref{eq:optimality_condition_h1_discrete} are enforced on vertices only.
This is due to the fact, that the expression $|T\vecv{u}_{h} - g_h|$ is not necessarily
cell-wise linear, even though $T\vecv{u}_{h} - g_h$ is.

Further we define: Let $F \in \Gamma$ be an oriented inner facet with adjacent cells $K_1, K_2 \in
\mathcal T$ and $\vec{\phi} \in L^2(K_1 \cup K_2)^{d\times m}$ with $\vec{\phi}|_{K_i} \in
C(K_i)^{d\times m}$, $i \in \{1, 2\}$, 
a square integrable function which allows continuous representations
$\vec{\phi}|_{K_i}$ on each cell $K_i$ individually.
We then define the jump term $[\vec{\phi}]_F \in C(F)^{d\times m}$ by
$[\vec{\phi}]_F(\vecv{x}) := \vec{\phi}|_{K_1}(\vecv{x}) - \vec{\phi}|_{K_2}(\vecv{x})$ and omit
the index as in $[\vec{\phi}]$ when the facet in question is clear.
For outer facets only the term corresponding to the existing adjacent cell is
considered.
Note that $[\vec{\phi}]_F$ is in general dependent on the orientation of $F$, while
e.g.\ $[\vec{n}^\transpose \vec{\phi}]_F$ for the oriented facet-normal $\vec{n}$ is not.

\subsection{On Image Interpolation Methods}\label{sec:ImageInterpolation}

For aligned grids as in \cref{fig:grid-discretization} image data $A \in
[0,1]^{n_1 \times n_2}$ is interpolated to a grid function $g_h \in Z_h$ by
nodal sampling.
This way, when resampling at pixel centers one receives the original image.
For non-aligned grids, 
i.e.\ when pixel centers do not correspond to vertices
of the grid, a mapping of the image $A$ to some grid function $g_h \in Z_h$ has
to be computed.
Choosing nodal sampling (e.g.\ of bilinearly interpolated image data) for this
map can result in aliasing and is a well known undesired effect in image
processing \cite{bredies2018mathematical}.
The Shannon-Whittaker sampling theorem provides ample conditions for correct
equidistant sampling of bandwidth-limited functions
\cite{shannon1948mathematical,bredies2018mathematical}, namely to remove
high-frequency contributions (e.g.\ by a linear smoothing filter) before
sampling on a coarser grid.
It is however unclear to us how to apply this to general unstructured sampling
points.

We evaluate three projection schemes and then explore one alternative.
Let $g\in C(\Omega)$ be the bilinear interpolation of the nodal image data
$A$.
For methods involving quadratures the quadrature over a cell $E \in \mathcal T$
is computed using a simple averaging quadrature on a Lagrange lattice of
degree $\lceil\op{diam}(E)\rceil$, i.e.\ we properly scale the number of
quadrature points depending on the size of the cell (and thus the number of
pixels it covers) to avoid aliasing effects.

The projection schemes are defined as follows:

\begin{enumerate}[(i)]
  \item
    \Code{nodal}:
    Nodal interpolation, i.e.\ $g_h \in Z_h$ such that $g_h(\vecv{x}) := g(\vecv{x})$ at every
    mesh vertex $\vecv{x}$.
  \item
    \Code{l2_lagrange}:
    Standard $L^2$-projection, i.e.\ $g_h \in Z_h$ such that $\int_\Omega g_h
    \cdot \phi \di[\vecv{x}] = \int_\Omega g \cdot \phi \di[\vecv{x}]$ for all $\phi \in
    Z_h$.
  \item
    \Code{qi_lagrange}:
    The general $L^1$-stable quasi-interpolation operator as proposed in
    \cite{ErnGuermond:17}, i.e.\ the continuous $g_h \in Z_h$ is given by
    setting its nodal degrees of freedom at each vertex to the arithmetic mean
    of the corresponding local degrees of freedom of a discontinuous interpolant
    $\widetilde{g}_h \in \widetilde{Z}_h := \{f \in L^\infty(\Omega) : f|_K \in P_1(K), K
      \in \mathcal T\}$, which is
    defined on each cell $K \in \mathcal T$ by its local nodal degrees of
    freedom $\sigma_{K,i} = \frac{1}{|K|} \int_K g(\vecv{x}) \cdot \rho_{K,i}(\vecv{x}) \di[\vecv{x}]$, $i
    = 1, \dotsc, 3$ where the test function $\rho_{K,i}$ in our case is given in
    barycentric coordinates $\lambda_K$ by
    $\rho_{K, i}(\lambda_K) = 12 \lambda_{K,i} - 3$.
    The reader may refer to \cite{ErnGuermond:17} for details on the general
    construction.
  \item
    \Code{l2_pixel}:
    Our proposed method, i.e.\ minimizing the sum of pointwise squared errors
    over all pixel coordinates:
    $\inf_{g_h \in Z_h} \sum_{\vecv{x} \in \Omega \cap \Z^2}
    \tfrac{1}{2} \nnew{|g_h(\vecv{x}) - g(\vecv{x})|}^2$.
\end{enumerate}

Note that \Code{l2_pixel} may be interpreted as an $L^2$-projection with
a cell-dependent averaging quadrature rule, adapted to the original image pixel
locations.

\begin{figure}
  \centering
  \includegraphics[width=0.2\textwidth]{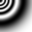}%
  \caption{Discrete input image for mesh interpolation comparison, $32 \times 32$ pixels.}
  \label{fig:mesh-interpolation-input}
\end{figure}

\begin{figure}
  \centering%
  \def\figurewidth{0.16\textwidth}%
  \foreach \meshsize in {32,16,13} {%
    \foreach \method in {nodal,l2_lagrange,qi_lagrange,l2_pixel} {%
      \begin{subfigure}[t]{\figurewidth}%
          \PlotImageMesh[width=\textwidth,height=\textwidth]%
          {image-mesh-interpolation\meshsize_\method.csv}%
      \end{subfigure}\hspace{0.1cm}%
    }\\%
  }%
  \caption{
    Interpolated finite element functions for varying mesh sizes and
    interpolation methods.
    Methods from left to right: \Code{nodal}, \Code{l2_lagrange}, \Code{qi_lagrange},
    \Code{l2_pixel}.
    Number of mesh vertices in each dimension from top to bottom: 32, 16, 13.
  }%
  \label{fig:mesh-interpolations}
\end{figure}

\begin{table}
\centering
\caption{PSNR values (higher is better) of interpolated mesh functions from
    \cref{fig:mesh-interpolations}, sampled at image coordinates. \nnew{Here \Code{nr_vertices} denotes the number of vertices per dimension.}}
\begin{tabular}{c|cccc}
\texttt{nr\_vertices} & \texttt{nodal}     & \texttt{l2\_lagrange} & \texttt{qi\_lagrange} & \texttt{l2\_pixel} \\ 
32                    & $\infty$    & 55.826       & 22.606 & 327.597          \\ 
16                    & 21.141      & 22.020       & 21.483 & 23.106 \\ 
13                    & 18.861      & 19.650        & 19.053 & 19.981          \\ 
\end{tabular}
\label{tab:mesh-interpolations-psnr}
\end{table}


\begin{table}[ht]
\centering
\caption{SSIM values (higher is better) of interpolated mesh functions from
    \cref{fig:mesh-interpolations}, sampled at image coordinates. \nnew{Here \Code{nr_vertices} denotes the number of vertices per dimension.}}
\begin{tabular}{c|cccc}
\texttt{nr\_vertices} & \texttt{nodal}     & \texttt{l2\_lagrange} & \texttt{qi\_lagrange} & \texttt{l2\_pixel} \\ 
32           & 1.00000   & 0.99998      & 0.97846      & 1.00000   \\ 
16           & 0.93868   & 0.95133      & 0.94757      & 0.96693   \\ 
13           & 0.89307   & 0.90888      & 0.90053      & 0.92080   \\ 
\end{tabular}\label{tab:mesh-interpolations-ssim}
\end{table}

We evaluate these schemes by interpolating the discrete image given in
\cref{fig:mesh-interpolation-input} onto $V_h$ for regular meshes of different
sizes.
In \cref{fig:mesh-interpolations} we see the interpolated results as cellwise
linear functions.
We sample the interpolated functions on the original image grid
$\Omega_h$ and quantify the error to the original image.
In particular,
\cref{tab:mesh-interpolations-psnr,tab:mesh-interpolations-ssim} list the $\op{PSNR}$, given for two scalar images
$u, v: \Omega_h \to [0,1]$ by
$\op{PSNR}(u, v):=-10 \log_{10}\big(\frac{1}{|\Omega_h|}\sum_{\vecv{x}\in\Omega_h}
(u(\vecv{x}) - v(\vecv{x}))^2 \big)$,
and the structural similarity index $\op{SSIM}$ as defined in~\cite{WaBoShSi:04}.
Perhaps unsurprisingly, due to construction, \Code{l2_pixel} produces
results closest to the original image for both metrics.
We note that, theoretically, \Code{l2_pixel} in the finest mesh setting should
provide exact results (i.e.\ infinite $\op{PSNR}$) and the discrepancy is due to
numerical imprecision.

It should be noted that while \Code{l2_pixel} seems to provide visually
superior results, it does not necessarily preserve other quantities of the
image, such as the total mass, which may be relevant for the image processing task
in question.
Further, it is not necessarily well-defined for meshes finer than the original
image and in that case a regularization needs to be applied.

We conclude that, since digital images are usually given as an array of
brightness values and the output of image processing algorithms is expected to
be so too, working in unstructured finite element spaces comes with the
drawback of information loss due to mesh interpolation.
Apart from increased complexity, this may be considered a problem for applying
unstructured adaptive finite element methods to image processing tasks in practice.

\subsection{Residual A-Posteriori Error Estimator}
\label{ssec:error-estimator-res}

In \cite{HintermullerRincon-Camacho:14} an a-posteriori error estimate for a
smooth functional, composed of an $L^2$-data term and a total
variation term, was derived using residual-based methods for variational
inequalities.
We try to use a similar technique for our general setting to establish a
suitable guiding criterion for adaptive spatial discretization. Our presentation avoids the variational inequality setting in
\cite{HintermullerRincon-Camacho:14} by making use of strong convexity, which seems to lead to a simpler presentation.
Additionally, we consider both choices of \cref{setting-identity,setting-nabla},
which requires different interpolation estimates depending on
the corresponding coercivity of $a_B$ stated in
\cref{assumption-a-coercivity}.


We recall
$E: V \to \overline\R$ from \eqref{eq:primal_smoothed} in the
split form $E(\vecv{u}) = F(\vecv{u}) + G(\Lambda \vecv{u})$, where
$\Lambda = (T, \nabla): V \to W$ as in \cref{ssec:predual-and-dualization} and
$F: V \to \overline\R$, $G: W \to \overline\R$ are given by
\begin{align*}
  F(\vecv{u}) &:=
    \tfrac{\alpha_2}{2} \|T\vecv{u} - g\|_{L^2}^2
    + \tfrac{\beta}{2} \|S\vecv{u}\|_{L^2}^2, \\
  G(\Lambda \vecv{u}) &:=
    \alpha_1 \int_\Omega \phi_{\gamma_1}(|T\vecv{u} - g|) \di[\vecv{x}]
    + \lambda \int_\Omega \phi_{\gamma_2}(|\nabla \vecv{u}|_F) \di[\vecv{x}].
\end{align*}
The optimality conditions in this setting, see \cite[Remark III.4.2]{EkelandTemam:99}, are given by
\begin{align*}
  -\Lambda^* \vec{p} &\in \partial F(\vecv{u}), \\
  \vec{p} &\in \partial G(\Lambda \vecv{u}),
\end{align*}
where we accounted for the sign change in the dual variable $\vec{p}$ to be
consistent with the formulation in \eqref{eq:optimality_condition_h1}.
We assume for the moment that $\alpha_1=0$ and the discrete solution pair $\vecv{u}_{h} \in V_h$, $\vec{p}_h \in W_h^*$
satisfies the
corresponding discrete optimality conditions
\begin{equation} \label{eq:opt_cond_discrete}
  \begin{aligned}
    -\Lambda^* \vec{p}_h &\in \partial F_h(\vecv{u}_{h}), \\
    \vec{p}_h &\in \partial G_h(\Lambda \vecv{u}_{h}),
  \end{aligned}
\end{equation}
where for clarity we explicitly denoted with $F_h: V_h \to \overline\R$, $G_h:
W_h \to \overline\R$ the restrictions of $F$, $G$ on the respective discrete spaces, since their subdifferentials
differ from those of $F$ and $G$ respectively. Note that $W_h$ denotes the (pre-)dual space of $W_h^*$.
One may obtain this system analogously to \cite{ours2022semismooth1} by application of
the Fenchel duality from \cite[Remark III.4.2]{EkelandTemam:99} on the corresponding discrete
energy functional. 
Indeed, \cref{thm:optcond_reg_nabla} in general allows for $V$ and $W$ to be
appropriate subspaces.
This way existence of a discrete solution pair $\vecv{u}_{h} \in V_h$, $\vec{p}_h \in W_h^*$
is ensured. 
Since $F$ is Fr\'echet-differentiable we see that
\begin{align*}
  \<\partial F(\vecv{u}), \vecv{v}\>_{V^*,V}
  &= \alpha_2 \<T\vecv{u} - g, T\vecv{v}\>_{L^2} + \beta \<S\vecv{u}, S\vecv{v}\>_{L^2} \\
  &= a_B(\vecv{u}, \vecv{v}) - l(\vecv{v})
\end{align*}
with $l(\vecv{v}) := \alpha_2\<T\vecv{v}, g\>_{L^2}$.
Analogously $\<\partial F_h(\vecv{u}_{h}), \vecv{v}_h\>_{V_h^*, V_h} = a_B(\vecv{u}_{h}, \vecv{v}_h) -
l(\vecv{v}_h)$ and from \eqref{eq:opt_cond_discrete} we infer
for all $\vecv{v}_h \in V_h$:
\begin{align} \label{eq:opt_cond_discrete_f}
  -\<\vec{p}_h, \Lambda \vecv{v}_h\>_{L^2} = a_B(\vecv{u}_{h}, \vecv{v}_h) - l(\vecv{v}_h).
\end{align}
For $\partial G$ on the other hand we have the following non-obvious
observation.

\begin{lemma} \label{lem:opt_cond_discrete_conformity}
  Let $\vec{p}_h = (p_{h,1}, \vec{p}_{h,2}) \in W_h^*$ satisfy
  \eqref{eq:opt_cond_discrete}.
  If $\alpha_1 = 0$ then
  \begin{align*}
    \vec{p}_{h} \in \partial G(\Lambda \vecv{u}_{h}),
  \end{align*}
  where $\vec{p}_{h} \in W_h^* \subset W^*$ is identified with an element of the
  dual space $(W^*)^* = W$.
  \begin{proof}
    Let $\vec{q} = (q_1, \vec{q}_2) \in \partial G(\Lambda \vecv{u}_{h})$.
    Due to \eqref{eq:optimality_condition_h1} and \eqref{eq:opt_cond_discrete} we
    deduce $q_1 = p_{h, 1} = 0$ since $\alpha_1 = 0$.
    From \eqref{eq:optimality_condition_h1} we infer for $\vec{q}_{2}$, and analogously for
    $\vec{p}_{h,2}$, due to \eqref{eq:opt_cond_discrete} that pointwise, whenever
    $\nabla \vecv{u}_{h} \neq 0$, we have
    \begin{align*}
      \vec{q}_2
      = \frac{\lambda \nabla \vecv{u}_{h}}{\max\{\gamma_2, |\nabla \vecv{u}_{h}|_F\}}
      = \vec{p}_{h,2}.
    \end{align*}
    Therefore $\vec{q}_{2}$ is piecewise constant whenever $\nabla \vecv{u}_{h} \neq 0$.
    Further both $\vec{q}_{2}$ and $\vec{p}_{h,2}$ are bounded pointwise by $|\vec{q}_{2}|_F \le \lambda$ and
    $|\vec{p}_{h,2}|_F \le \lambda$.
    On cells where $\nabla \vecv{u}_{h} = 0$ holds, $\vec{q}_{2}$ may assume any value bounded
    by $|\vec{q}_{2}|_F \le \lambda$.
    Since this holds true for $\vec{p}_{h,2}$ in particular, one may choose $\vec{q} :=
    (0, \vec{p}_{h,2}) \in \partial G(\Lambda \vecv{u}_{h})$.
  \end{proof}
\end{lemma}

We note that for $\alpha_1>0$ the optimality conditions \eqref{eq:opt_cond_discrete} might not be valid for $\vec{p}_{h} \in W_h^*$ due to the presence of the operator $T$. A possibility to circumvent this shortcoming might be to enlarge the space $W_h^*$ to $\{(p_1, \vec{p}_2) \in L^2(\Omega) \times L^2(\Omega)^{d \times m}:
        p_1|_K \in P_1(K),
        \vec{p}_2|_{K} \in P_0(K)^{d \times m}, K \in \mathcal T \}$. 
However, even in such a new setting \cref{lem:opt_cond_discrete_conformity} would not necessarily hold
true for $\alpha_1 > 0$, as
$\partial G(\Lambda \vecv{u}_{h}) = \alpha_1 \tfrac{T \vecv{u}_{h} - g}{\max\{\gamma_1, |T
\vecv{u}_{h} - g|\}}$ does not need to be piecewise linear, even if $\vecv{u}_{h}$ is.
We will continue the derivation nevertheless, noting that theoretical
justification is lacking for $\alpha_1 > 0$. That is in the sequel we assume that 
\begin{equation}\label{eq:opt_cond_discrete_conformity}
\vec{p}_{h} \in \partial G(\Lambda \vecv{u}_{h}).
\end{equation}
  
With these observations we are ready to estimate the error \nnew{$\|\vecv{u}_{h} - \vecv{u}\|_V$.}

\nnew{
\begin{theorem}\label{Thm:ErrorEstimate}
Let $\vecv{u} \in V$ be a solution of \eqref{eq:primal_smoothed}, and let $\vecv{u}_{h} \in V_h$ and $\vec{p}_{h}\in W_h^*$ satisfy \eqref{eq:opt_cond_discrete}. 
If \cref{assumption-a-coercivity} and \eqref{eq:opt_cond_discrete_conformity} hold, then we have
\begin{enumerate}[(i)]
\item for \cref{setting-nabla}
\begin{align*}
  \|\vecv{u}_{h} - \vecv{u}\|_{H^1(\Omega)^m}^2 &\le C \bigg(\sum_{K \in \mathcal T} h_K^2 \|\alpha_2 T^*(T\vecv{u}_{h} - g) +
    T^*p_{h,1}
    - \beta \Laplace \vecv{u}_{h} - \Div \vec{p}_{h,2}\|_{L^2(K)^m}^2 \\
  &\qquad\qquad + \sum_{F\in\Gamma} h_F \|[\vec{n}^\transpose(\beta\nabla \vecv{u}_{h} + \vec{p}_{h,2})]\|_{L^2(F)^m}^2 \bigg)
\end{align*}
with constant $C > 0$, and
\item for \Cref{setting-identity}
\begin{align*}
  \|\vecv{u}_{h} - \vecv{u}\|_{L^2(\Omega)^m}^2
  &\le C \bigg( \sum_{K \in \mathcal T} \|\alpha_2 T^*(T\vecv{u}_{h} - g) + T^* p_{h,1} +
    \beta \vecv{u}_{h} - \Div \vec{p}_{h,2}\|_{L^2(K)^m}^2  \\
  &\qquad\qquad + \sum_{F\in\Gamma} h_F^{-1} \|[\vec{n}^\transpose\vec{p}_{h,2}]\|_{L^2(F)^m}^2 \bigg)
\end{align*}
for some constant $C > 0$.
\end{enumerate}
\end{theorem}}
\nnew{The proof of this theorem is postponed to Appendix \ref{Sec:ProofErrorEstimate}.
}

\subsubsection*{Error Indicators}
\nnew{\Cref{Thm:ErrorEstimate} gives raise to a local error indicator. In particular,} for a cell $K \in \mathcal T$ we define our local error indicators as follows
\begin{equation} \label{eq:errind}
\eta_{\nnew{\mathrm{res}},K}^2 := \tilde \eta_{\nnew{\mathrm{res}},K}^2
  + \sum_{\substack{F \in \Gamma \\ F \cap K \neq \emptyset}} \tilde \eta_{\nnew{\mathrm{res}},F}^2,
\end{equation}
where $\tilde \eta_{\nnew{\mathrm{res}},K}$ and $\tilde \eta_{\nnew{\mathrm{res}},F}$ are chosen as follows.
For \cref{setting-nabla} and $S=\nabla$ we set 
\begin{align} \label{eq:errind_nabla}
  \tilde \eta_{\nnew{\mathrm{res}},K}^2 &:= h_K^2 \big\|\alpha_2 T^*(T\vecv{u}_{h} - g) + T^*p_{h,1} - \beta \Laplace \vecv{u}_{h} -
  \Div \vec{p}_{h,2} \big\|^2_{L^2(K)^m}, \\
  \tilde \eta_{\nnew{\mathrm{res}},F}^2 &:= h_{F}  \big\|[ \vec{n}^\transpose(\beta \nabla \vecv{u}_{h} + \vec{p}_{h,2}) ]
  \big\|^2_{L^2(F)^m}, \nonumber
\end{align}
where for piecewise linear $\vecv{u}_{h}$ and piecewise constant $\vec{p}_{h,2}$ the
terms $\Laplace \vecv{u}_{h}$ and $\Div \vec{p}_{h,2}$ in \eqref{eq:errind_nabla} vanish.
For \cref{setting-identity} and $S=I$ we set
\begin{align*}
  \tilde \eta_{\nnew{\mathrm{res}},K}^2 &:= \big\|\alpha_2 T^*(T\vecv{u}_{h} - g) + T^* p_{h,1} + \beta \vecv{u}_{h} - \Div
  \vec{p}_{h,2}\big\|_{L^2(K)^m}^2, \\
  \tilde \eta_{\nnew{\mathrm{res}},F}^2 &:= h_{F}^{-1} \big\|[\vec{n}^\transpose \vec{p}_{h,2}]\big\|_{L^2(F)^m}^2.
\end{align*}
For \cref{setting-identity} with $S = I$ the facet indicator $\tilde \eta_{\nnew{\mathrm{res}},F}$ scales inversely with the
diameter and is therefore not very useful in the context of adaptive refinement.
This is due to $a_B$ only being $L^2$-coercive 
 in that case, which
limits our choice of interpolation estimates and it is
unclear to us whether this result can be improved.
Showing efficiency of these estimators may be considered in future work and is
expected to work out in a similar way as in
\cite{HintermullerRincon-Camacho:14}.

In the sequel we call the indicator \eqref{eq:errind} \textit{residual based indicator}.

\subsection{Adaptive Refinement Strategy}\label{sec:refinement}
For refinement we bisect triangles using the newest-vertex strategy \cite{NoSiVe:09}, i.e.\ the
bisection edge is chosen to be opposite of the vertex which was inserted last.
In an adaptive refinement setting we mark cells for refinement using the
greedy D\"orfler marking strategy \cite{NoSiVe:09} with \nnew{$\theta_{\text{mark}}\in [0,1]$, e.g.\ $\theta_{\text{mark}} = 0.5$},
i.e.\ given error indicators in descending order $(\eta_{K_n})_{1 \le n \le
|\mathcal T|}$ for triangles $K_n \in \mathcal T$ we refine the first
$n_{\text{mark}} \in \N$ triangles that satisfy
\begin{align*}
  \sum_{n = 1}^{n_{\text{mark}}} \eta_{K_n} \ge \theta_{\text{mark}}
  \sum_{n=1}^{|\mathcal T|} \eta_{K_n}.
\end{align*}

Our generic adaptive finite element algorithm for solving \eqref{eq:primal_smoothed} is presented in \cref{alg:AFEM}. 

\begin{algorithm}\caption{AFEM for \eqref{eq:primal_smoothed}} \label{alg:AFEM}~\\[-0.5cm]
\begin{enumerate}[1.)]
\item Choose an initial grid $\mathcal T$.\label{point1}
\item Find a finite element solution $\vecv{u}_{h}$ of \eqref{eq:primal_smoothed} on $\mathcal T$, e.g.\ by \cite[Algorithm 1]{ours2022semismooth1}.\label{point2}
\item \nnew{Compute on each cell a local indicator, e.g.\ given by 
\eqref{eq:errind} or \eqref{eq:errind-pd}}.\label{point3}
\item Stop if a certain termination criterion holds (e.g.\ 
a certain number of refinements is reached%
), or continue with next step. \label{point4}
\item Determine $\mathcal R_h$, the set of cells to be refined, by using the greedy D\"orfler strategy as described in \cref{sec:refinement}.\label{point5}
\item Refine the cells $K\in \mathcal R_h$ by utilizing newest-vertex bisection and set  $\mathcal{T}$ to be the new grid.\label{point6}
\item Project the image data onto $\mathcal{T}$ (e.g.\ by \textnormal{\Code{l2_lagrange}} or \textnormal{\Code{l2_pixel}}) and continue with step 2.).\label{point7}
\end{enumerate}

\end{algorithm}

\section{Applications}\label{sec:numerical-examples}

In the following, we demonstrate that our model \eqref{eq:primal_smoothed}, in conjunction with \cref{alg:AFEM}, can effectively be applied in practice to solve image processing tasks such as image inpainting and optical flow computation in image sequences.
Thereby we apply \cref{alg:AFEM} in the following way: 
As initial grid in step~\ref{point1} a coarse grid is chosen which is coarser than the observed, original image grid. Based on the local indicators and the greedy D\"orfler strategy in step~\ref{point6} only a certain amount of, namely $n_{\text{mark}}$, cells are refined. Utilizing \cref{alg:AFEM} in this way generates a grid which is partly coarse and partly fine, see e.g.\ \cref{fig:inpainting-mesh}. In this vein we call this strategy an adaptive coarse-to-fine scheme. Moreover, in all our experiments to find a finite element solution $\vec{u}_h$ in step~\ref{point2} of \cref{alg:AFEM} we utilize the primal-dual semi-smooth Newton method of \cite{ours2022semismooth1}. 

The implementation of the performed algorithms is done in the programming
language Julia \cite{BeEdKaSh:17} and can be found at \cite{HilbAFEMCode:23}.
All evaluations were done on a notebook with a Intel Core i7-10875H CPU @ 2.30GHz and 32GB RAM.

\subsection{Inpainting} \label{ssec:examples-inpainting}
\begin{figure}
  \centering
  \includegraphics[width=0.3\textwidth]{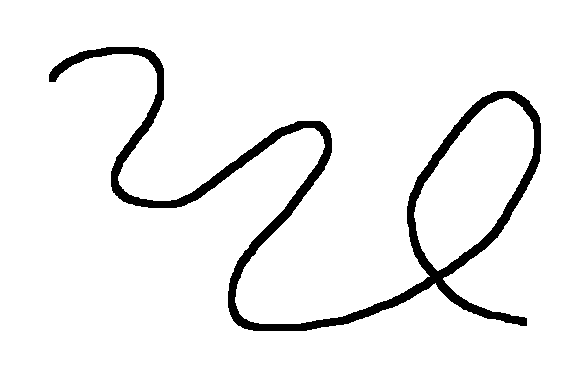}
  \includegraphics[width=0.3\textwidth]{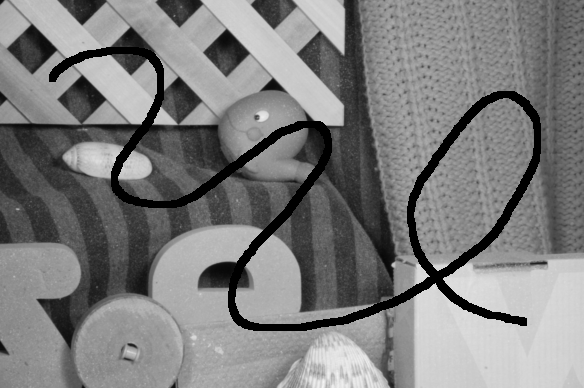}
  \caption{%
    Inpainting domain (left), corrupted image (right).}
  \label{fig:inpainting}
\end{figure}
\begin{figure}
  \centering
  \includegraphics[width=0.3\textwidth]{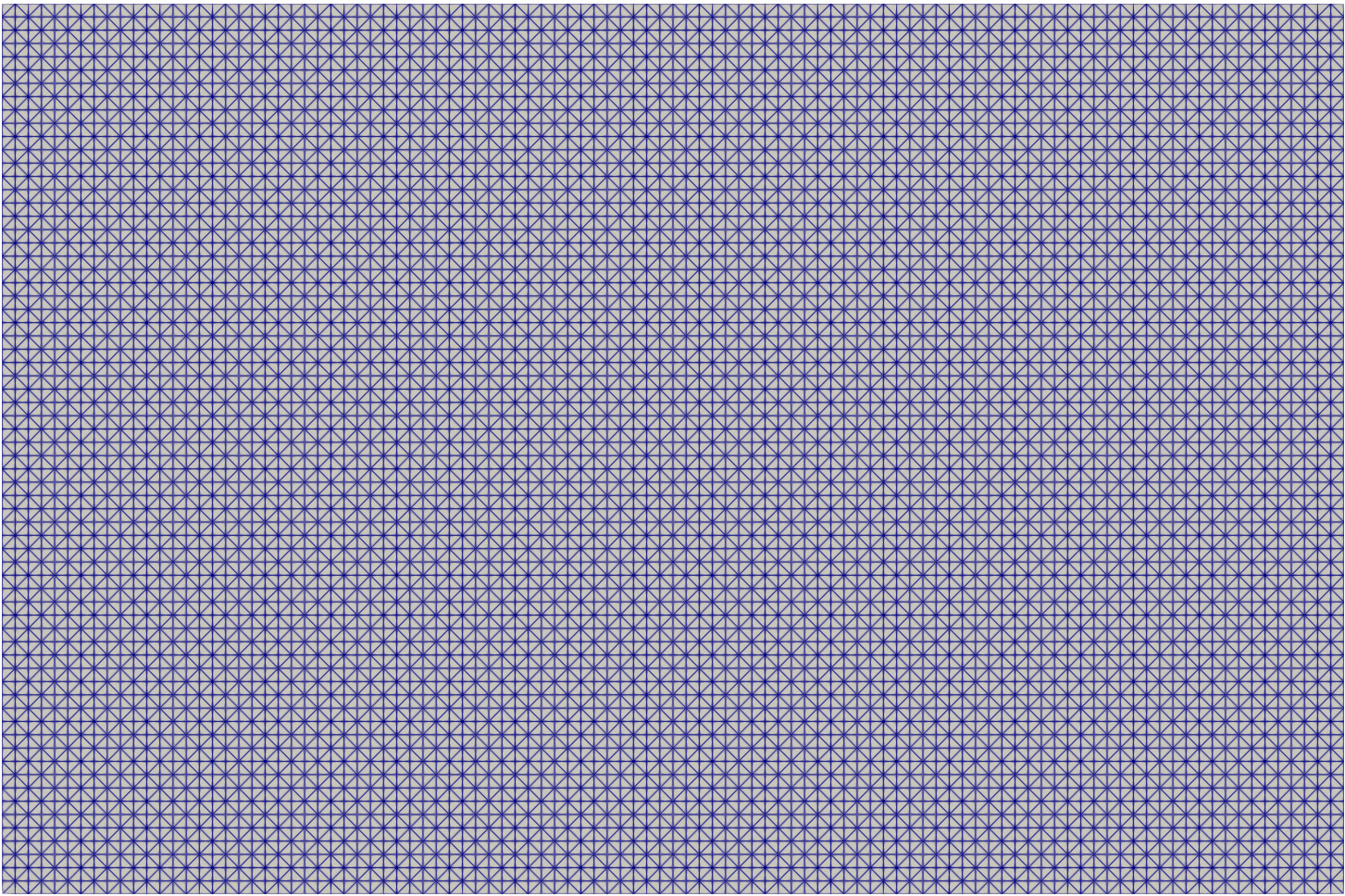}
  \includegraphics[width=0.3\textwidth]{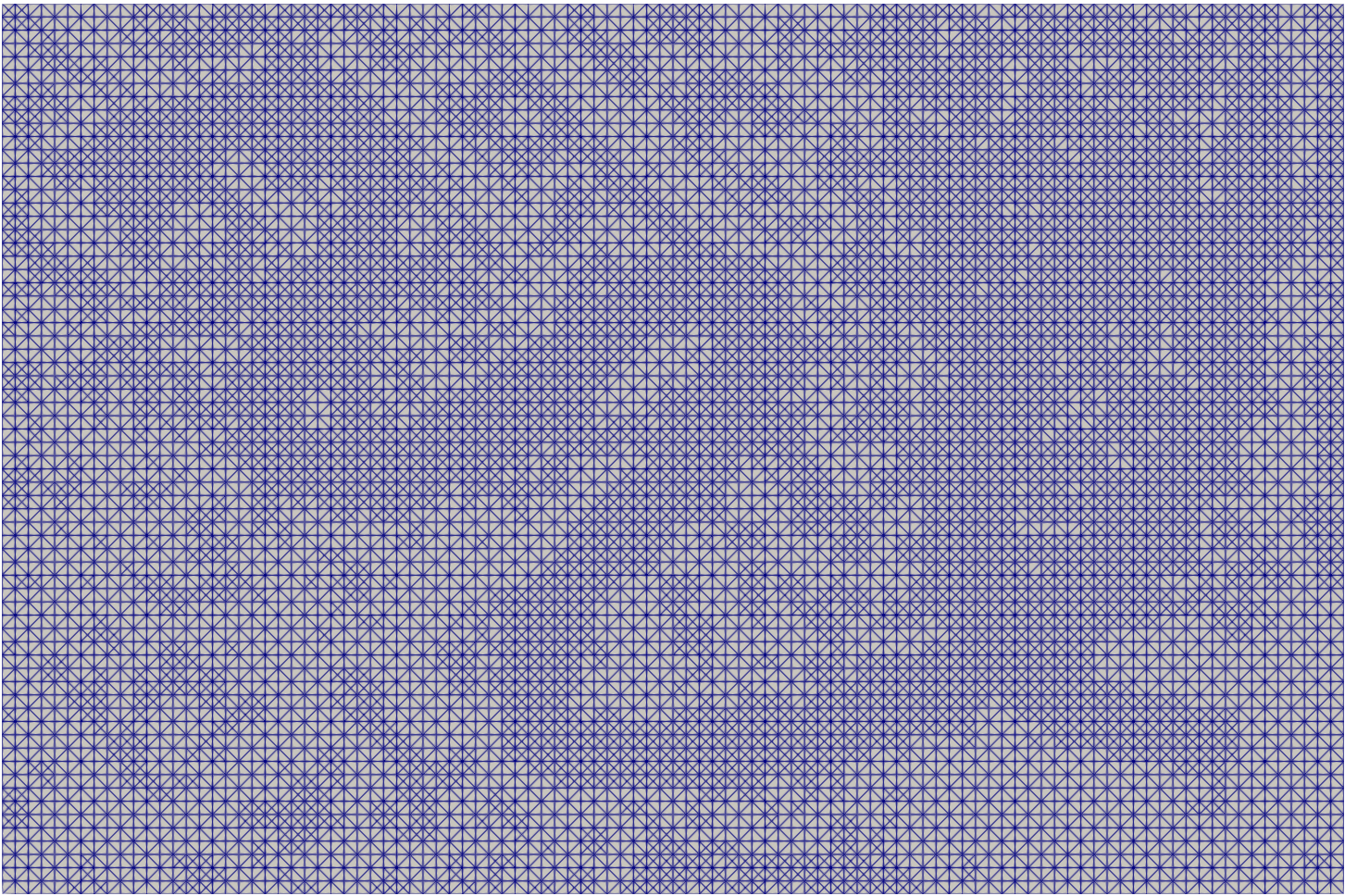}
  \includegraphics[width=0.3\textwidth]{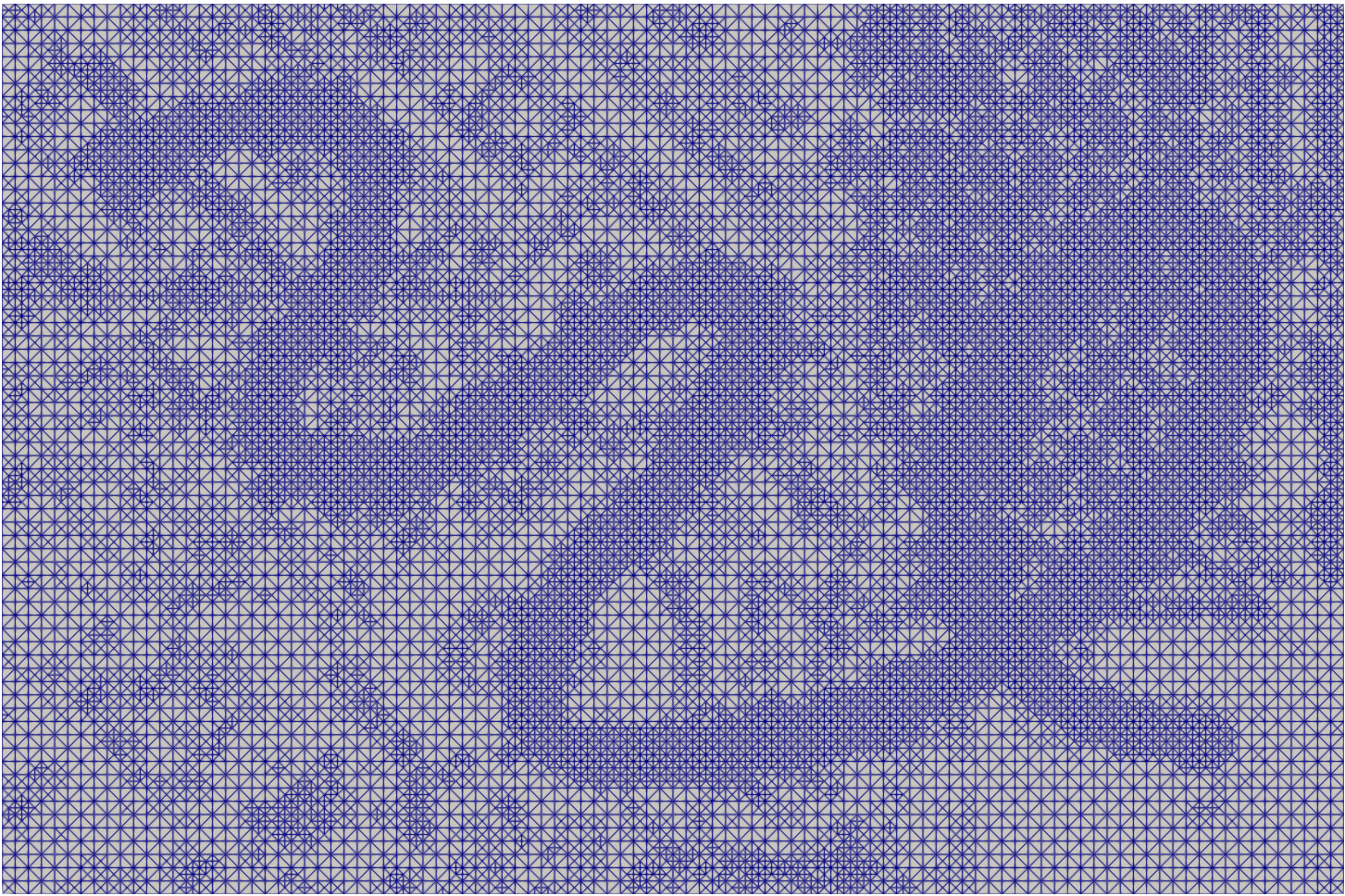}
  \includegraphics[width=0.3\textwidth]{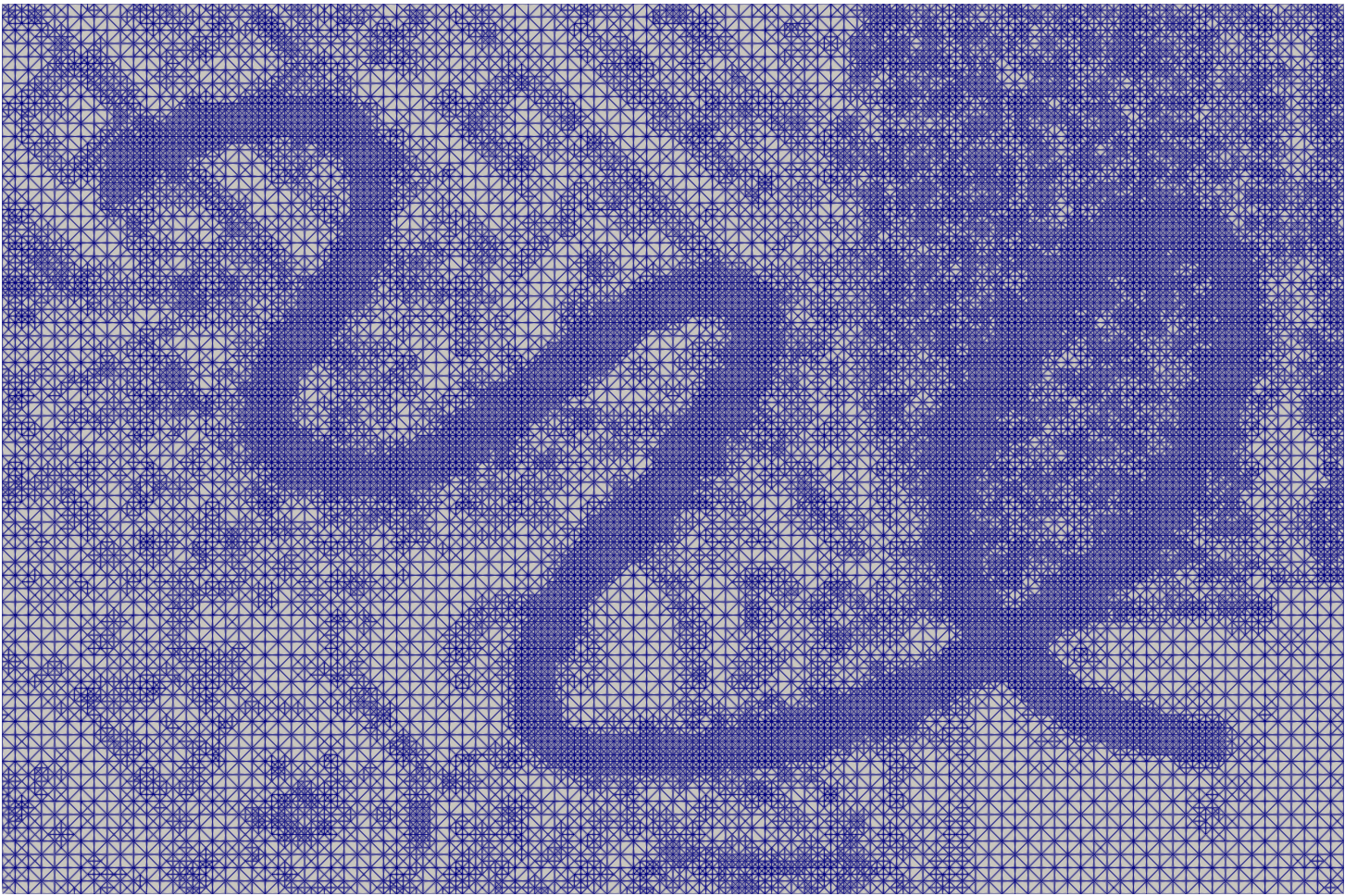}
  \includegraphics[width=0.3\textwidth]{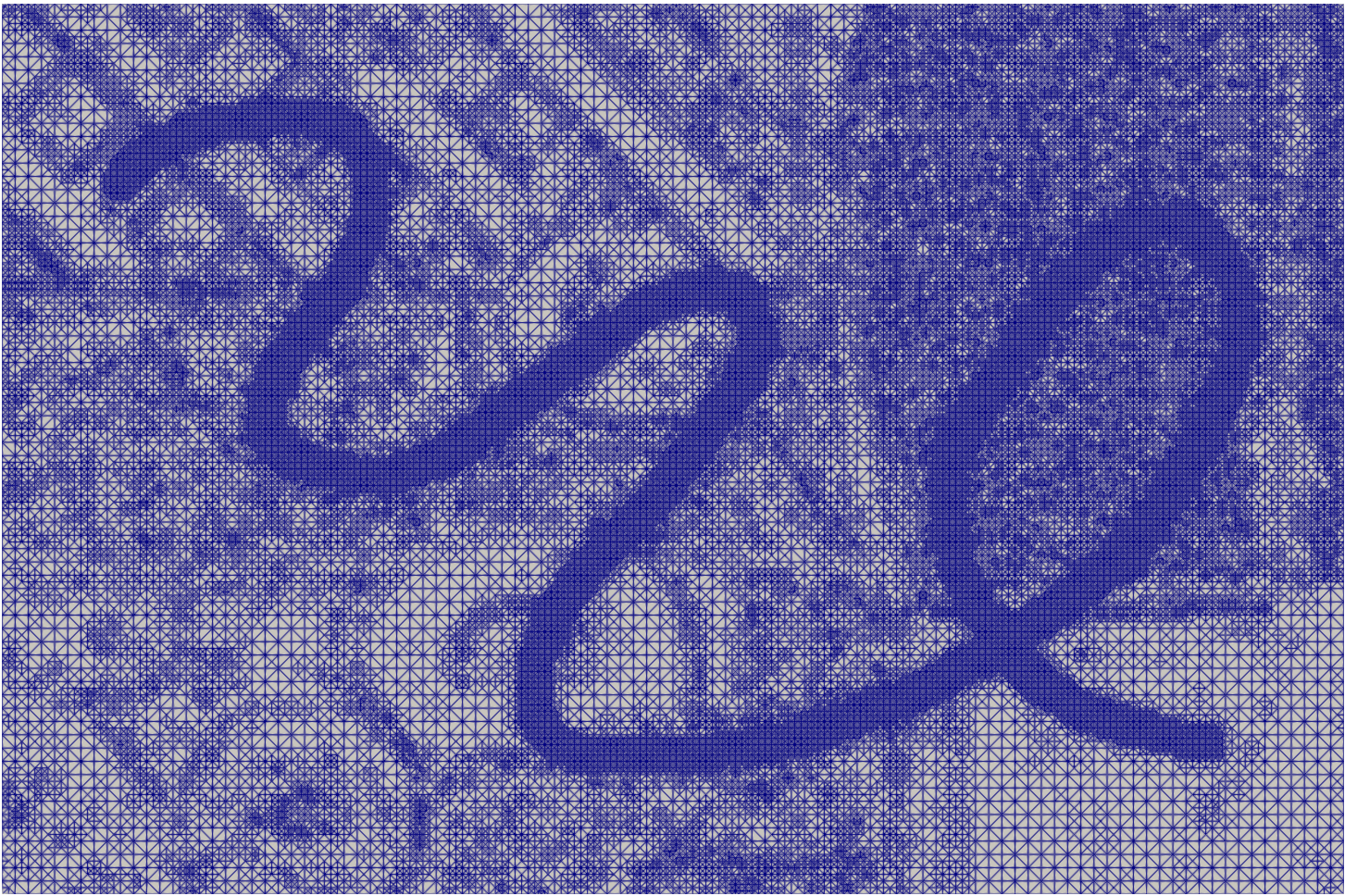}
  \includegraphics[width=0.3\textwidth]{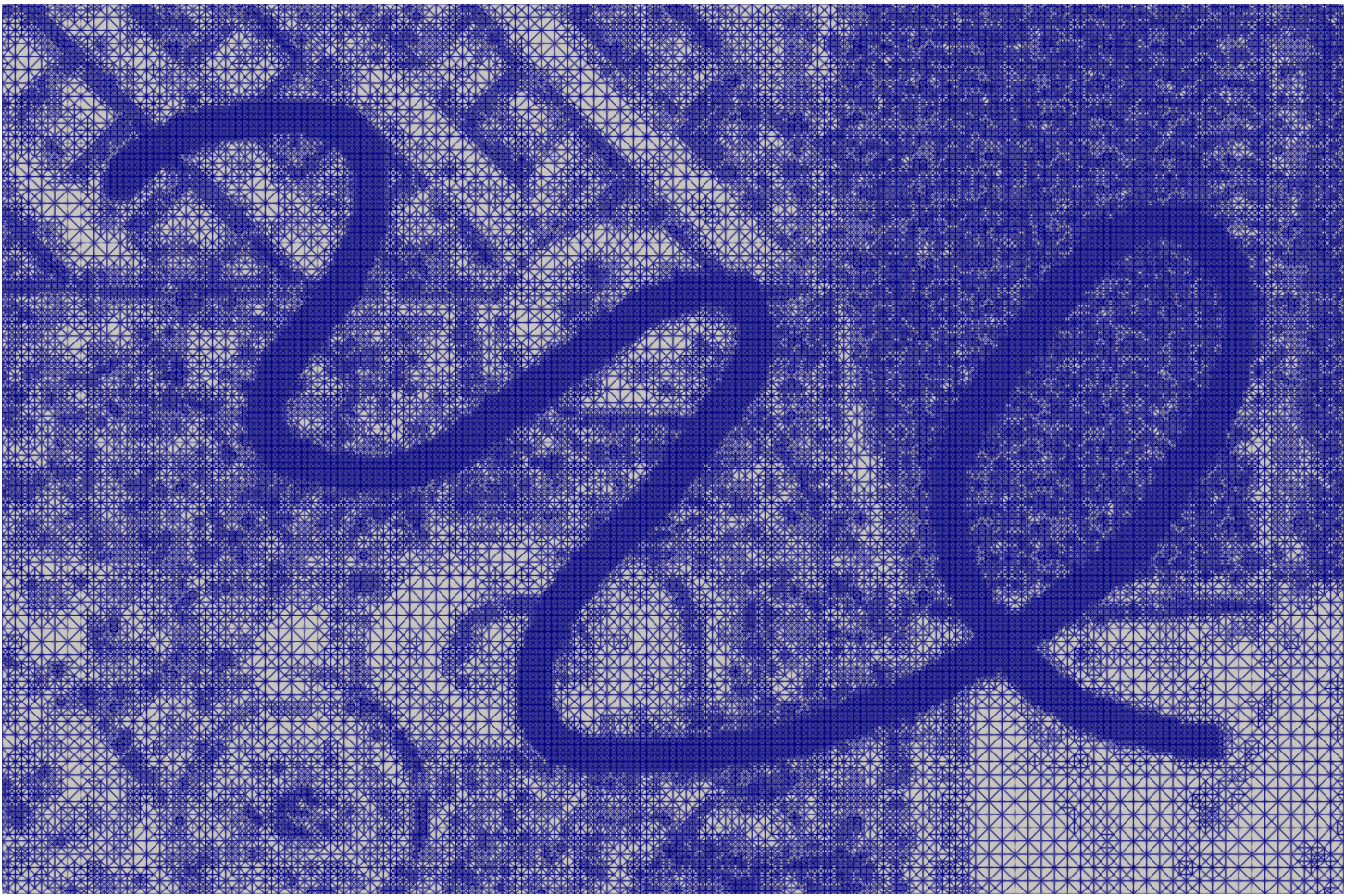}
  \caption{%
    Adapted mesh for inpainting during the iteration for $n_{\text{coarsen}}=5=n_{\text{refine}}$ \nnew{using the residual based indicators and $\theta_{\text{mark}} = 0.5$}.
  }
  \label{fig:inpainting-mesh}
\end{figure}

Image inpainting is the problem of reconstructing a given defected image $g \in L^2(\Omega\setminus D)$, whereby any image information in $D \subset \Omega$ is missing. The domain $D$ is therefore called inpainting domain. \Cref{fig:inpainting} depicts an example of an inpainting domain $D$ and a respective corrupted observation $g$. We tackle this inpainting task by using model \eqref{eq:primal_smoothed} and utilizing \cref{alg:AFEM}. Thereby we set $T:=\Id_{\Omega\setminus D}$, where
$$
(\Id_{\Omega\setminus D} u)(\vec{x}):=\begin{cases}
u(\vec x) & \text{if } \vec{x}\in \Omega\setminus D,\\
0& \text{if } \vec{x}\in D,
\end{cases}
$$
for any $u\in L^2(\Omega)$, $S=I$ and we manually select the parameters $\alpha_1 = 0$, 
$\alpha_2 = 50$, 
$\lambda = 1$,
$\beta = 1\cdot 10^{-5}$,
$\gamma_1 = 1 \cdot 10^{-4}$,
$\gamma_2 = 1 \cdot 10^{-4}$.
The finite element solution in each iteration is computed by \cite[Algorithm 1]{ours2022semismooth1} with 
$\epsilon_{\text{newton}} =  1 \cdot 10^{-4}$.
The initial grid is chosen to be uniform and $\tfrac{1}{2^{n_{\text{coarsen}}/2}}$ 
of the image resolution, rounded down to integer values, where $n_{\text{coarsen}}\in\N_0$. The algorithm performs then $n_{\text{refine}} = \lfloor\log_2\tfrac{\text{Number of cells in the original image}}{\text{Number of cells in the initial grid}}\rfloor$ refinements, see \cref{fig:inpainting-mesh} for
$n_{\text{coarsen}}= 5$
 which leads to 
$n_{\text{refine}}= 5$.
This allows that the smallest element is approximately of the size of the elements of the original image, but not smaller. 
\nnew{Note that our indicators, presented in \cref{sec:finite-elements}, may not be insightful in the inpainting domain. For example in the case of the primal-dual based indicators if $n_{\text{coarsen}}$ is so large that the initial grid consists only of the inpainting domain, then the local indicators calculate to 0 and no element will be marked for refinement. To circumvent such situations, after applying the greedy Dörfler marking strategy we additionally mark all elements in the inpainting area for refinement for both types of indicators.}
We use \Code{qi_lagrange} for projecting onto the finite element space except when $n_{\text{refine}}=0$, in which case we use nodal interpolation of the image data. 
The choice of \Code{qi_lagrange} over \Code{l2_pixel} and \Code{l2_lagrange} is based on the following factors. 
Firstly, \Code{l2_pixel} is not well-suited for inpainting tasks because it is a global operator, which can cause data leakage from the inpainting domain. 
Additionally, \Code{qi_lagrange} seems to perform slightly better than \Code{l2_lagrange} in our example, offering faster results and improved restoration quality.

\begin{figure}
  \centering
  \includegraphics[width=0.3\textwidth]{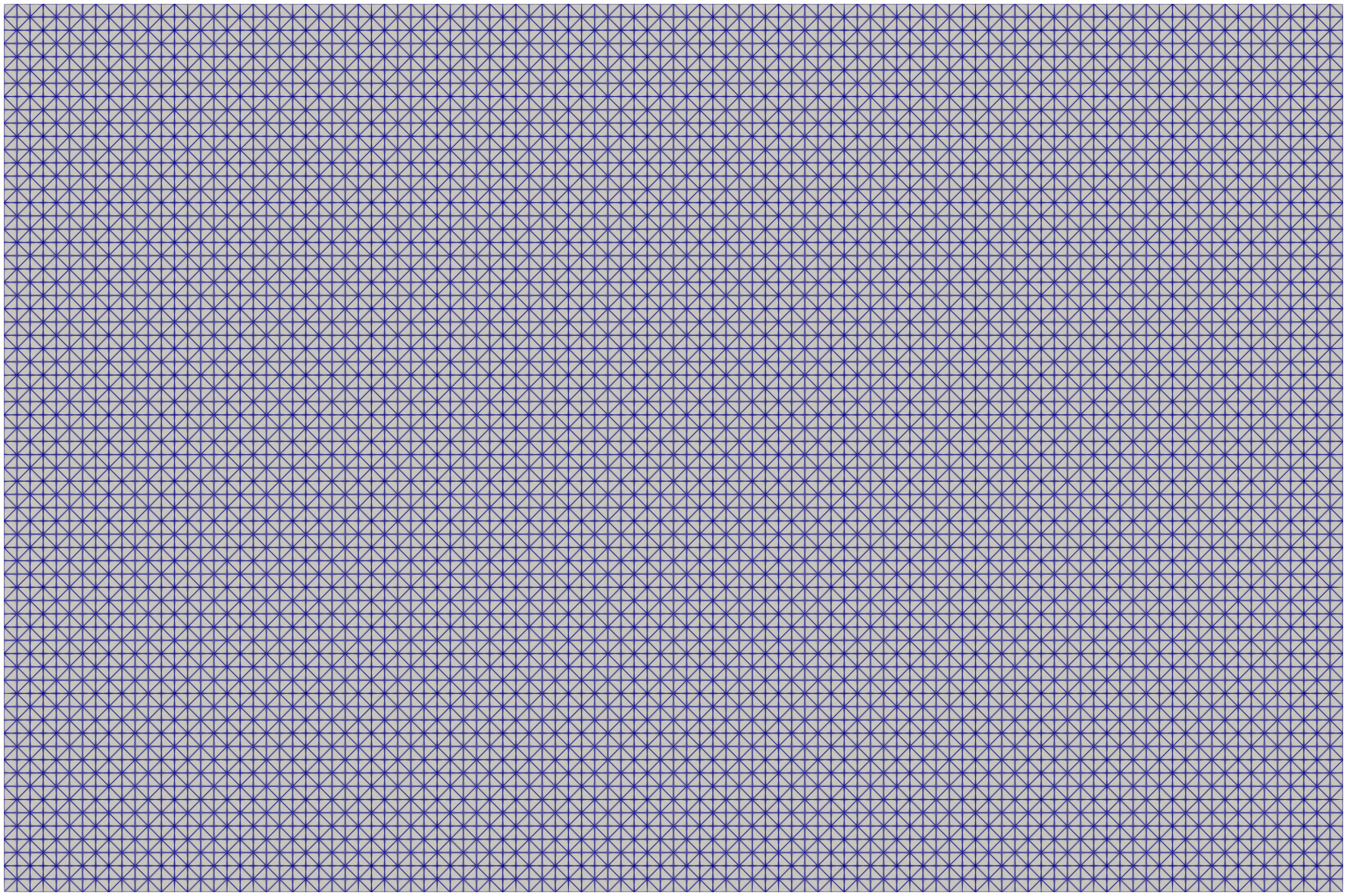}
  \includegraphics[width=0.3\textwidth]{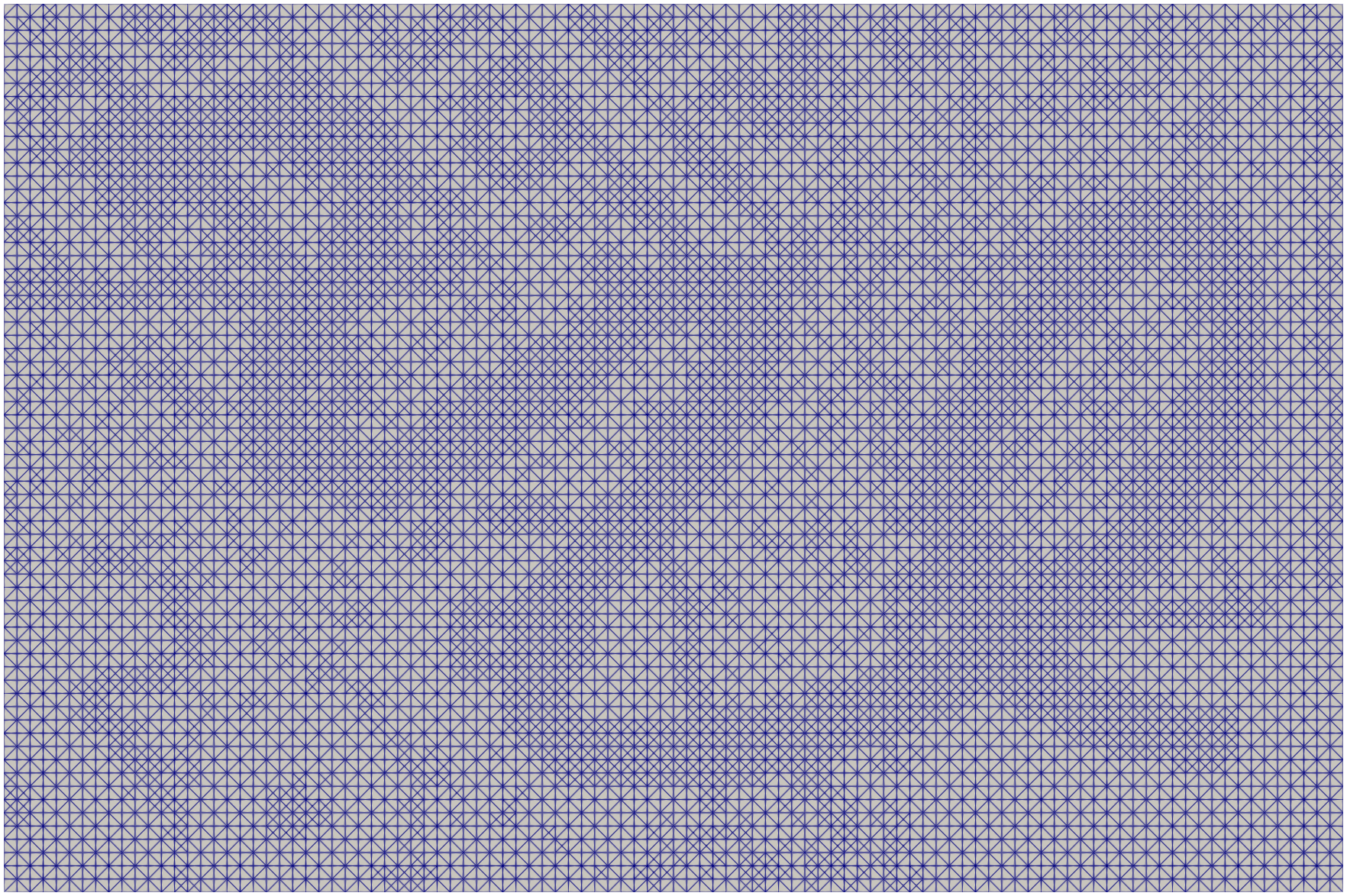}
  \includegraphics[width=0.3\textwidth]{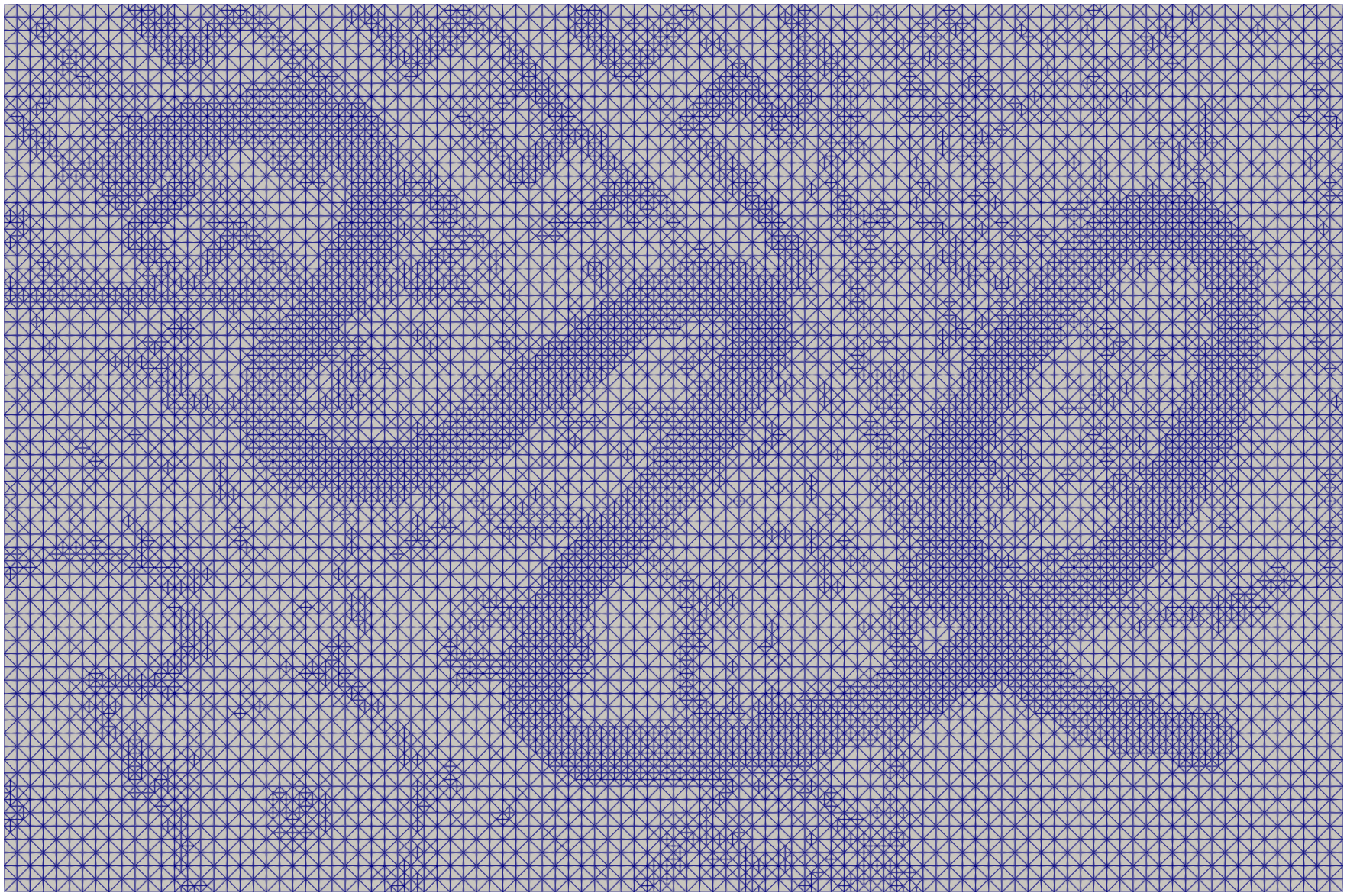}
  \includegraphics[width=0.3\textwidth]{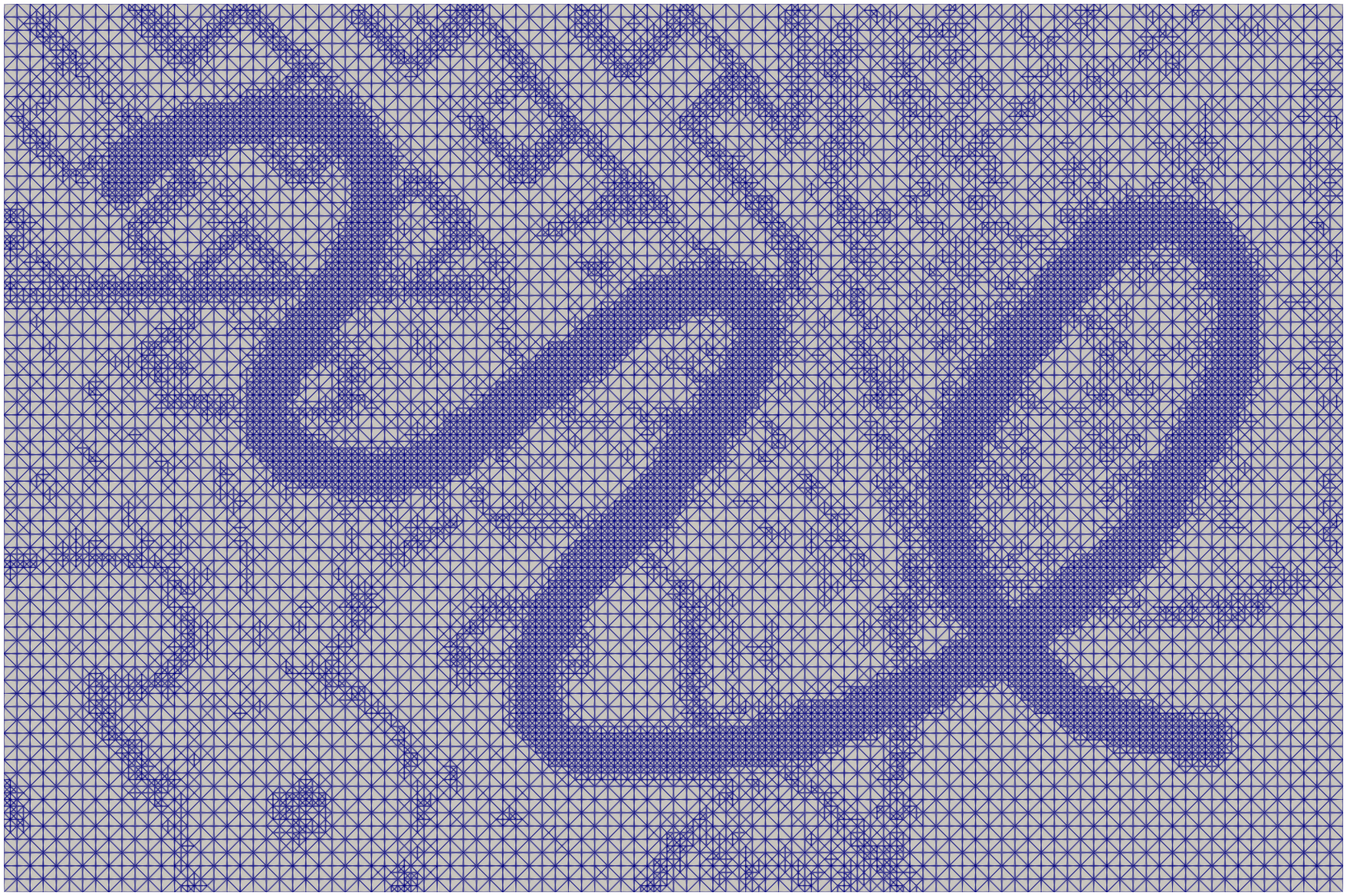}
  \includegraphics[width=0.3\textwidth]{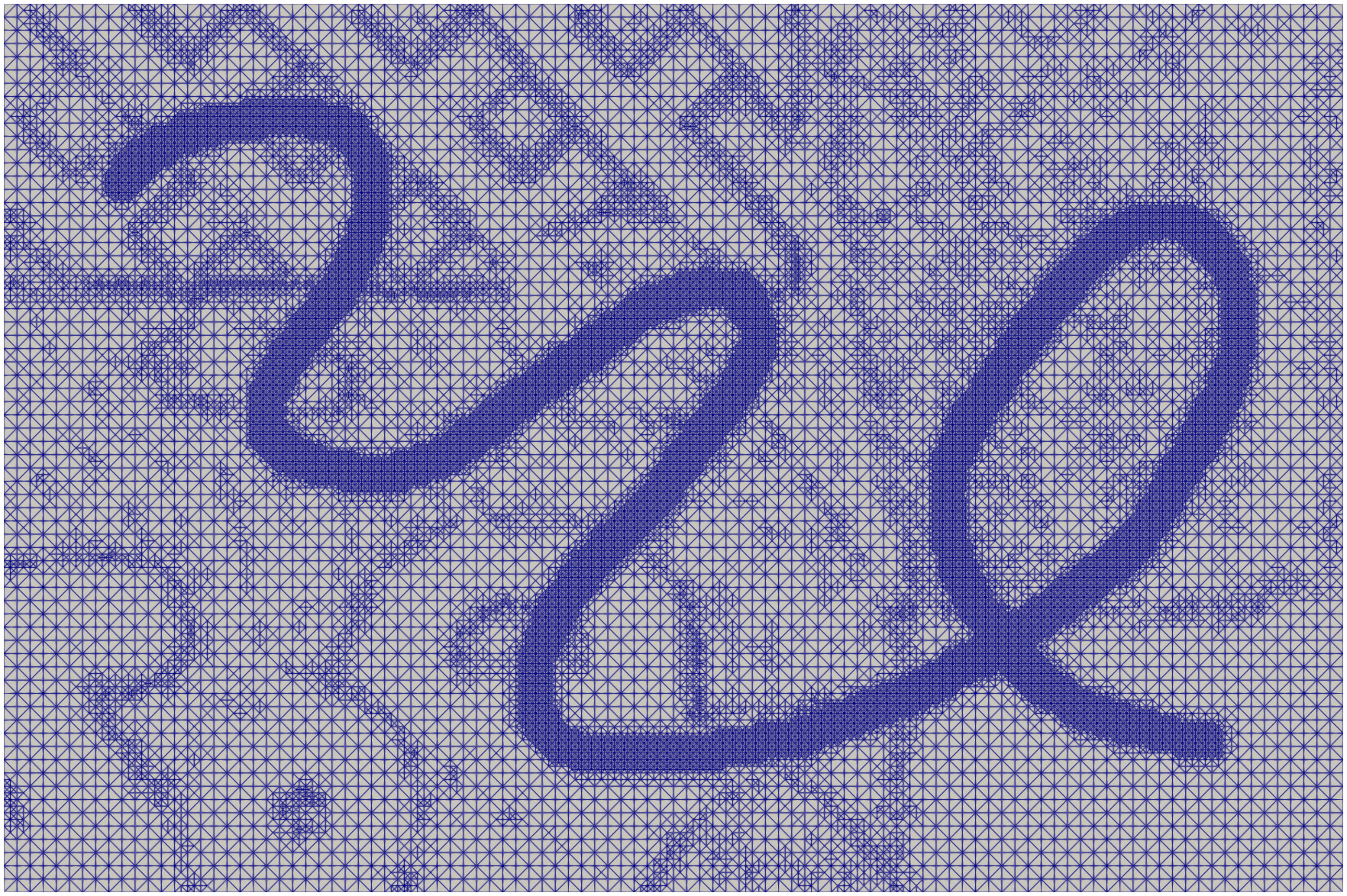}
  \includegraphics[width=0.3\textwidth]{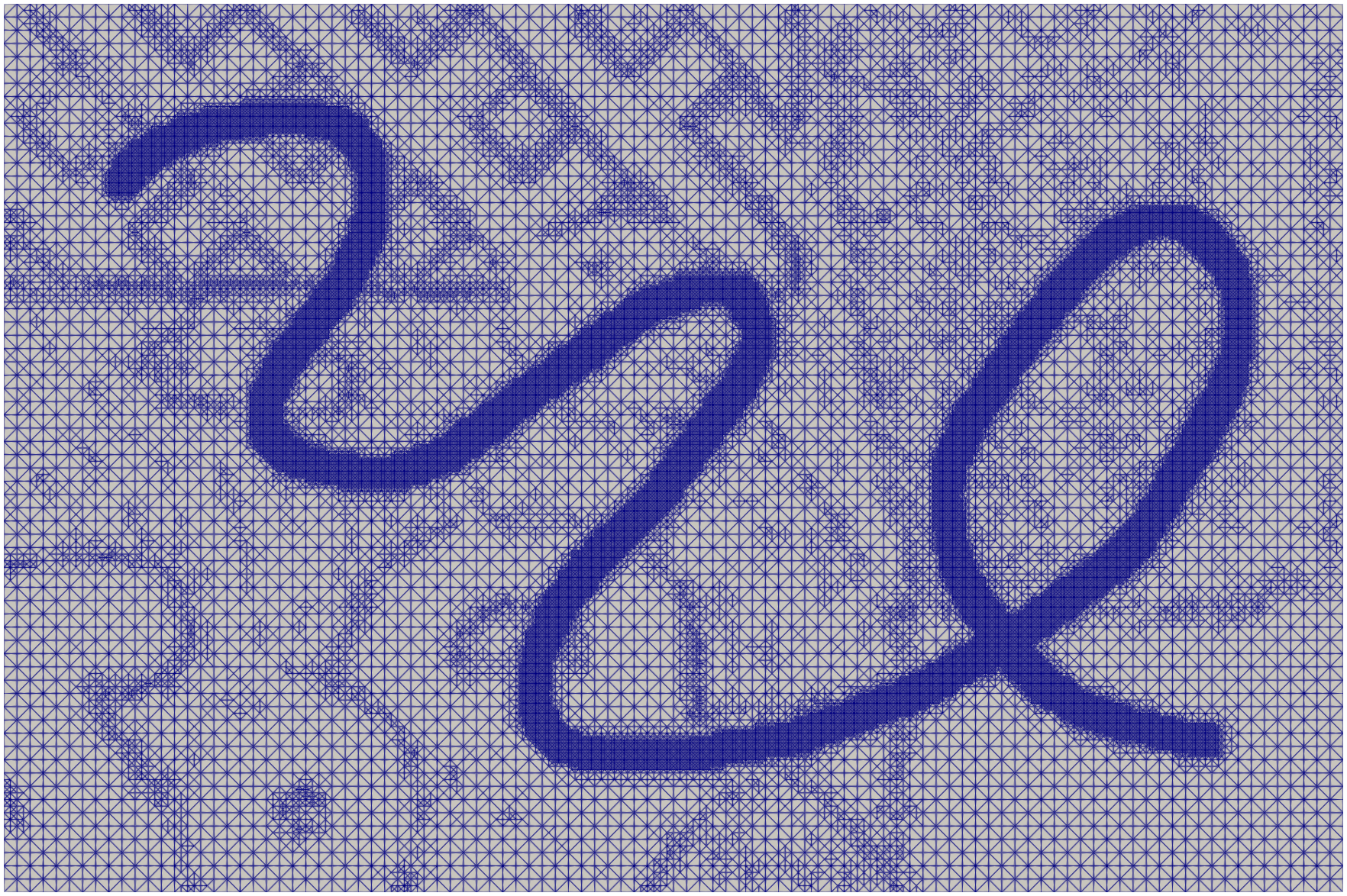}
  \caption{%
    \nnew{Adapted mesh for inpainting during the iteration for $n_{\text{coarsen}}=5=n_{\text{refine}}$ with the primal-dual based indicators and $\theta_{\text{mark}} = 0.99$.}
  }
  \label{fig:inpainting-mesh-pd}
\end{figure}
\begin{figure}
  \centering
  \includegraphics[width=0.3\textwidth]{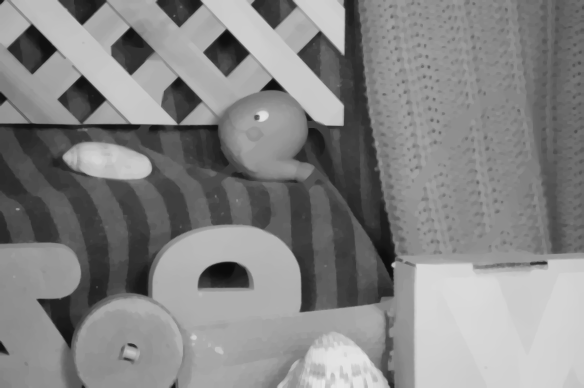}
  \includegraphics[width=0.3\textwidth]{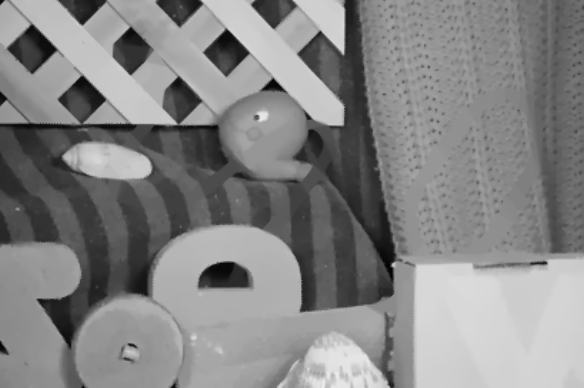}
  \includegraphics[width=0.3\textwidth]{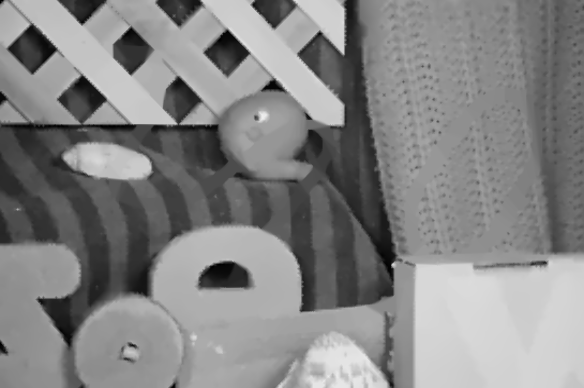}
  \caption{%
  Reconstruction obtained with $n_{\text{coarsen}} = 0$ (left), $n_{\text{coarsen}} = 5$ (middle), $n_{\text{coarsen}} = 10$ (right) \nnew{using the residual based indicators}.
  }
  \label{fig:inpainting:reconstruction}
\end{figure}
\begin{figure}
  \centering
  \includegraphics[width=0.3\textwidth]{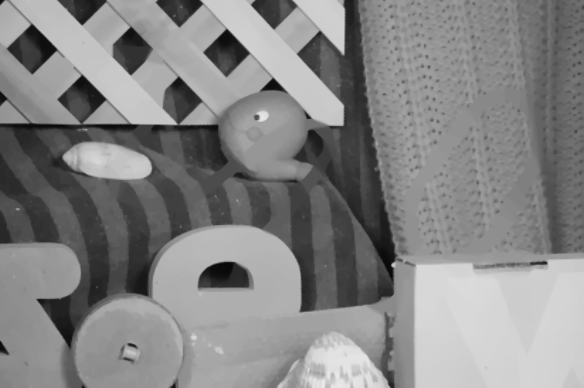}
  \includegraphics[width=0.3\textwidth]{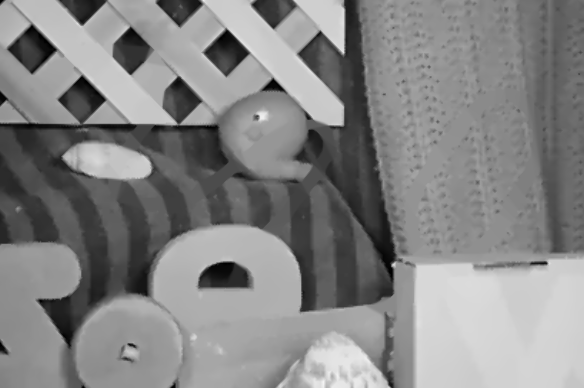}
  \includegraphics[width=0.3\textwidth]{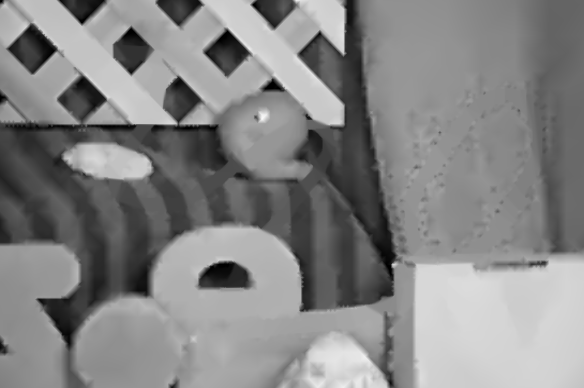}
  \caption{%
  \nnew{Reconstruction obtained with $n_{\text{coarsen}} = 1$ (left), $n_{\text{coarsen}} = 5$ (middle), $n_{\text{coarsen}} = 10$ (right) using the primal-dual based indicators.}
  }
  \label{fig:inpainting:reconstruction-pd}
\end{figure}

\nnew{We present in \cref{fig:inpainting-mesh} for the residual based indicators with $\theta_{\text{mark}} = 0.5$ and in \cref{fig:inpainting-mesh-pd} for the primal-dual based indicators with $\theta_{\text{mark}} = 0.99$ (see below for the reason of this choice)} the adaptive mesh generated by \cref{alg:AFEM} in each iteration for
$n_{\text{refine}}= 5$. 
 Since the edges of the elements are depicted in blue in \cref{fig:inpainting-mesh,fig:inpainting-mesh-pd}, dark areas indicate fine refinement, while in brighter areas the mesh is coarser. \nnew{Recall that all elements in the inpainting domain are always marked for refinement}. Moreover, we observe that in uniform areas as in the right lower corner of our image, the mesh stays coarse as this seems sufficient to keep the approximation error at a reasonable level and still obtain a visually appealing reconstruction, see \cref{fig:inpainting:reconstruction,fig:inpainting:reconstruction-pd}.  

\nnew{Note that using the primal-dual gap as an error estimator by computing it on a discrete mesh requires solving the model in the primal-dual gap equation on a finer mesh than what the candidate solution is using in order to provide meaningful results. Otherwise on the same mesh it will evaluate to zero.
Moreover, we observed that for this application and with $\theta_{\text{mark}}=0.5$ the primal-dual based indicators allow to select for refinement only cells in the inpainting domain leading to undesirable results. 
However, setting $\theta_{\text{mark}} = 0.99$ leads to similar results as with the  residual based indicators and $\theta_{\text{mark}} = 0.5$, see \cref{fig:inpainting:reconstruction,fig:inpainting:reconstruction-pd}. 
From the mesh depicted in \cref{fig:inpainting-mesh-pd} one easily sees that even with $\theta_{\text{mark}} = 0.99$ the algorithm mainly refines in the inpainting domain. 
Hence the mesh is much coarser (76082 cells in the final grid with $n_{\text{coarsen}} = 5$) than the one generated with the help of the residual based indicators (152081 cells in the final grid with $n_{\text{coarsen}} = 5$). 
Nevertheless, the computation with the primal-dual based indicators is much more expensive, as the indicators need to be computed on a globally refined mesh in each iteration. 
In particular, for $n_{\text{coarsen}} = 5$ it takes $95.35$ seconds generating a restoration with $\op{PSNR} = 17.07$ and $\op{SSIM} = 0.7459$. 
\pgfkeys{/pgf/number format/.cd,fixed,fixed zerofill,precision=2}
On the contrary, using the residual based indicators for $n_{\text{coarsen}} = 5$ a result is generated within $57.01$ seconds with $\op{PSNR} = 17.12$ and $\op{SSIM} = 0.7968$. Quantitative similar results are obtained with different $n_{\text{coarsen}}$, suggesting that for inpainting our residual based indicators perform better in \cref{alg:AFEM} than our primal-dual based indicators.
}


\subsubsection*{Influence on the initial grid}

\begin{table}
\centering
\caption{Comparison of different values $n_{\text{coarsen}}$ for solving the inpainting problem in \cref{fig:inpainting} using the residual based indicators.}\label{tab:inpainting}
\begin{tabular}{c c|c c c c }
\toprule
$n_{\mathrm{coarsen}}$ & $n_{\mathrm{refine}}$ & {\# cells in final grid} & {PSNR}  & {SSIM}  & {Time (s)} \\ \hline \hline
0  & 0  & \numprint{451242} & 17.20 & 0.8083 & 65.05 \\ 
1  & 1  & \numprint{319120} & 17.14 & 0.8184 & 83.41 \\ 
2  & 2  & \numprint{253678} & 17.13 & 0.8132 & 81.33 \\ 
3  & 3  & \numprint{207341} & 17.14 & 0.8081 & 72.76 \\ 
4  & 4  & \numprint{177321} & 17.12 & 0.8027 & 64.93 \\ 
5  & 5  & \numprint{152081} & 17.12 & 0.7968 & 57.01 \\ 
6  & 6  & \numprint{135825} & 17.10 & 0.7908 & 50.30 \\ 
7  & 7  & \numprint{119846} & 17.08 & 0.7836 & 49.35 \\ 
8  & 8  & \numprint{109463} & 17.08 & 0.7758 & 42.09 \\ 
9  & 9  & \numprint{99839}  & 17.06 & 0.7659 & 36.25 \\ 
10 & 10 & \numprint{93638}  & 17.05 & 0.7598 & 33.29 \\ 
11 & 11 & \numprint{80006}  & 17.03 & 0.7499 & 29.22 \\ 
12 & 12 & \numprint{82627}  & 17.04 & 0.7537 & 30.26 \\ 
13 & 13 & \numprint{67377}  & 17.02 & 0.7385 & 35.15 \\ 
14 & 15 & \numprint{95524}  & 17.06 & 0.7603 & 35.31 \\ 
15 & 16 & \numprint{70524}  & 17.02 & 0.7434 & 26.39 \\ 
\bottomrule
\end{tabular}
\end{table}

\nnew{Here we restrict ourselves to the residual based indicators, since in inpainting they are superior to the primal-dual based indicators.}
We investigate the behavior of our algorithm for different choices of  $n_{\text{coarsen}}$, namely $n_{\text{coarsen}}\in \{0,1,2, \ldots,15\}$, when applied to the inpainting problem in \cref{fig:inpainting}. 
Note that with $n_{\text{coarsen}} =0$ no coarsening is performed and hence also no refinement, i.e.\ with this choice we indicate the performance of the algorithm on a uniform mesh as presented in \cite{ours2022semismooth1}. 
Since the size of the image in \cref{fig:inpainting} is $584 \times 388$ pixels, $n_{\text{coarsen}} = 15$ is the maximal possible number for coarsening the image yielding an image \nnew{with $3\times 2$ vertices}, i.e.\ 4 cells. 
Our findings, namely the number of cells in the finest grid, the PSNR and SSIM values of the restoration as well as the CPU time needed to obtain the reconstruction, are summarized in 
\cref{tab:inpainting} and visualized in \cref{fig:inpainting:investigation}. 
As expected, we observe that increasing $n_{\text{coarsen}}$ leads to a significant reduction in the number of cells in the finest, final grid. At the same time, the PSNR and SSIM values of the reconstruction remain relatively stable, resulting in a similarly appealing outcome, see \cref{fig:inpainting:reconstruction}. Nevertheless, from \cref{tab:inpainting} we observe that the quality measures slightly decrease with increasing $n_{\text{coarsen}}$, indicating a minor reduction in restoration capability. This might be attributed to the coarser grid, which may not allow the restored image to be represented with the same level of precision.
Considering the CPU times in \cref{tab:inpainting} we would like to emphasize that for $n_{\text{coarsen}}\geq 1$ a significant amount of time is spent on the refinement process. In particular, the more cells are refined the slower the process. This is the reason why there is a jump between $n_{\text{coarsen}}=0$ and $n_{\text{coarsen}}= 1$ in \cref{fig:inpainting:investigation}(b). Hence the reduction in CPU time with increasing $n_{\text{coarsen}}$ is then two-sided. First, due to the reduction of cells the used semi-smooth Newton method converges faster and second, the refinement process needs to refine less cells and hence is also faster. 
This results in the proposed adaptive algorithm terminating faster as $n_{\text{coarsen}}$ increases, compared to when $n_{\text{coarsen}}= 0$, albeit with a minor reduction in restoration ability. 

\tikzexternaldisable
\begin{figure}
\centering
\begin{subfigure}[t]{0.45\columnwidth}
\centering
\begin{tikzpicture}[scale=0.8]
\begin{axis}[
    xlabel={$n_{\mathrm{coarsen}}$},
    grid=both,
    legend style={at={(1.05,1)},anchor=north west},
    width=6.6cm, height=5cm,
    xmax=15.5
]
\addplot coordinates {
    (0,451242) (1,319120) (2,253678) (3,207341) (4,177321)
    (5,152081) (6,135825) (7,119846) (8,109463) (9,99839)
    (10,93638) (11,80006) (12,82627) (13,67377) (14,95524) (15,70524)
};
\end{axis}
\end{tikzpicture}
\caption{Number of cells in the final grid}
\end{subfigure}\hfill
\begin{subfigure}[t]{0.45\columnwidth}
\centering
\begin{tikzpicture}[scale=0.8]
\begin{axis}[
    xlabel={$n_{\mathrm{coarsen}}$},
    grid=both,
    legend style={at={(1.05,1)},anchor=north west},
    width=6.6cm, height=5cm,
    xmax=15.5
]
\addplot coordinates {
    (0,65.05) (1,83.41) (2,81.33) (3,72.76) (4,64.93)
    (5,57.01) (6,50.30) (7,49.35) (8,42.09) (9,36.25)
    (10,33.29) (11,29.22) (12,30.26) (13,35.15) (14,35.31) (15,26.39)
};
\end{axis}
\end{tikzpicture}
\caption{CPU time in seconds}
\end{subfigure}
\begin{subfigure}[t]{0.45\columnwidth}
\centering
\begin{tikzpicture}[scale=0.8]
\begin{axis}[
    xlabel={$n_{\mathrm{coarsen}}$},
    grid=both,
    legend style={at={(1.05,1)},anchor=north west},
    width=6.5cm, height=5cm,
    ymin=0, ymax=30,
    xmax=15.5
]
\addplot coordinates {
    (0,17.20) (1,17.14) (2,17.13) (3,17.14) (4,17.12)
    (5,17.12) (6,17.10) (7,17.08) (8,17.08) (9,17.06)
    (10,17.05) (11,17.03) (12,17.04) (13,17.02) (14,17.06) (15,17.02)
};
\end{axis}
\end{tikzpicture}
\caption{PSNR}
\end{subfigure}\hfill
\begin{subfigure}[t]{0.45\columnwidth}
\centering
\begin{tikzpicture}[scale=0.8]
\begin{axis}[
    xlabel={$n_{\mathrm{coarsen}}$},
    grid=both,
    legend style={at={(1.05,1)},anchor=north west},
    width=6.5cm, height=5cm,
    ymin=0, ymax=1,
    xmax=15.5
]
\addplot coordinates {
    (0,0.8083) (1,0.8184) (2,0.8132) (3,0.8081) (4,0.8027)
    (5,0.7968) (6,0.7908) (7,0.7836) (8,0.7758) (9,0.7659)
    (10,0.7598) (11,0.7499) (12,0.7537) (13,0.7385) (14,0.7603) (15,0.7434)
};
\end{axis}
\end{tikzpicture}
\caption{SSIM}
\end{subfigure}

\caption{Number of cells, CPU time, PSNR and SSIM for different initial image resolutions coarsened by the factor $\tfrac{1}{2^{n_{\text{coarsen}}/2}}$, $n_{\text{coarsen}}\in\{0,1,2,\ldots,15\}$.}
\label{fig:inpainting:investigation}
\end{figure}
\tikzexternalenable

\subsection{Motion Estimation}
 
Motion estimation is the problem of numerically computing the apparent motion field in a sequence of images. Following the discussion in \cite{ours2022semismooth1} the motion field $\vecv{u}$ of two (smooth) images $f_0$, $f_1$ may be computed solving model \eqref{eq:primal_smoothed} with 
\begin{align*}
  T \vecv{u} := \nabla f_{w} \cdot \vecv{u}, \qquad  g := \nabla f_{w} \cdot \vecv{u}_{0} - (f_{w} - f_0),
\end{align*}
where $f_w(\vecv{x}) := f_1(\vecv{x} + \vecv{u}_{0}(\vecv{x}))$ is a warped version of $f_1$ with $\vecv{u}_{0}$ being some smooth initial guess and $\vecv{x}\in\Omega$. To account for large displacements an optical flow warping algorithm is introduced \cite[Algorithm 3]{ours2022semismooth1}. 
We combine the warping technique from \cite[Algorithm 3]{ours2022semismooth1} with adaptive refinement leading to \cref{alg:optflow_adapt}.
\begin{algorithm}\caption{Optical flow algorithm with adaptive warping}~
    \label{alg:optflow_adapt}\\
  \Parameters{warping threshold $\epsilon_{\text{warp}}$, parameters for \cite[Algorithm 1]{ours2022semismooth1}}\\
  \Input{images $f_0, f_1$, initial guess $\vecv{u}_0$}\\
  \Output{motion fields $(\vecv{u}_k)_k$}

  \begin{algorithmic}
    \State{$f_{w,0}(\vecv{x}) = f_1(\vecv{x} + \vecv{u}_{0}(\vecv{x}))$}
    \For{$k = 1, 2, \dotsc$}
    \State{find approximate discrete solution $\vecv{u}_k$ to
        \eqref{eq:primal_smoothed} using \cite[Algorithm 1]{ours2022semismooth1}
        }
      \State{$f_{w,k}(\vecv{x}) = f_1(\vecv{x} + \vecv{u}_{k}(\vecv{x}))$}
      \If{$\tfrac{\|f_{w,k-1} - f_0\|_{L^2} - \|f_{w,k} -
      f_0\|_{L^2}}{\|f_{w,k-1} - f_0\|_{L^2}} < \epsilon_{\text{warp}}$}
        \State{refine mesh and reproject image data}
      \EndIf
    \EndFor
  \end{algorithmic}
\end{algorithm}

To approximately solve for $\vecv{u}_k$ in \cref{alg:optflow_adapt}, we use our
model \eqref{eq:primal_smoothed} by choosing $T \vecv{u}:=\nabla f_{w,k-1} \cdot \vecv{u}$ and $g := \nabla f_{w,k-1} \cdot \vecv{u}_{k-1}
-(f_{w,k-1} - f_0)$ and utilize the primal-dual semi-smooth Newton method of \cite{ours2022semismooth1}. Thereby we use the discrete space $Z_h$ for the images $f_0$, $f_1$, $f_{w,k}$, except for the warping step, in which $f_{w,k}(\vecv{x}) = f_1(\vecv{x} + \vecv{u}_k(\vecv{x}))$ is carried out at original image resolution by evaluating $f_1$ using bicubic interpolation and in a second step projected (e.g.\ using \Code{l2_lagrange} or \Code{l2_pixel}) onto the current finite element space in order to capture
more detailed displacement information. For $g$ we instead use a cellwise linear
discontinuous space to capture the discontinuous component $\nabla f_{w,k-1}$.

Loosely speaking, \cref{alg:optflow_adapt} solves the linearized optical flow equation for $u_k$
and warps the input data by the computed flow field until it no longer improves
on the data difference $f_{w,k} - f_0$, which indicates displacement.
In that case, the mesh is refined using the indicators from \eqref{eq:errind}
\nnew{or \eqref{eq:errind-pd} together with $\theta_{\text{mark}}=0.5$ in the greedy Dörfler marking strategy} and the process repeats, now including more detailed image data.

We note, that this approach to adaptivity allows us to start off with a coarse
mesh and refine cells only if deemed necessary by the error indicator.
In that respect it is different from the only other adaptive finite element
methods for
optical flow we are aware of, see \cite{belhachmi2011control,BelHec2016}, where the
mesh is initialized at fine image resolution first and iteratively coarsened only
after a costly computation of the flow field and a suitable metric for
adaptivity on this fine mesh has been established.

\subsubsection*{Experiments}


\begin{figure}
  \def\imgfraction{0.1428}
  \centering
  \foreach \data in {Dimetrodon,Grove2,Grove3,Hydrangea,RubberWhale,Urban2,Urban3,Venus} {
    \begin{subfigure}[t]{\imgfraction\columnwidth}
      \includegraphics[width=\textwidth]{fofm\data f0.png}
    \end{subfigure}%
    \begin{subfigure}[t]{\imgfraction\columnwidth}
      \includegraphics[width=\textwidth]{fofm\data f1.png}
    \end{subfigure}%
    \begin{subfigure}[t]{\imgfraction\columnwidth}
      \includegraphics[width=\textwidth]{fofm\data g.png}
    \end{subfigure}%
    \begin{subfigure}[t]{\imgfraction\columnwidth}
      \includegraphics[width=\textwidth]{fofm_w_c_adaptivevanilla\data output.png}
    \end{subfigure}%
    \begin{subfigure}[t]{\imgfraction\columnwidth}
      \includegraphics[width=\textwidth]{fofm_w_c_adaptivewarping\data output.png}
    \end{subfigure}%
    \begin{subfigure}[t]{\imgfraction\columnwidth}
      \includegraphics[width=\textwidth]{fofm_w_c_adaptiveadaptive-warping\data output.png}
    \end{subfigure}%
    \begin{subfigure}[t]{\imgfraction\columnwidth}
      \includegraphics[width=\textwidth]{fofm\data ground_truth.png}
    \end{subfigure}%
  }
  \caption{%
    Middlebury Optical Flow Benchmark for $S=\nabla$ using the residual based indicators: columns from left to right: $f_0$,
    $f_1$, image difference $f_1 - f_0$, computed
    optical flow $\vecv{u}$ using \cite[Algorithm 3]{ours2022semismooth1} without warping and with warping, computed optical flow $\vecv{u}$ using \cref{alg:optflow_adapt} (adaptive warping), ground truth optical flow . 
    Benchmarks from top to bottom: Dimetrodon, Grove2, Grove3, Hydrangea,
    RubberWhale, Urban2, Urban3, Venus.
  }
  \label{fig:middlebury}
\end{figure}

\begin{sidewaystable}
\scriptsize
  \centering
\caption{%
    Quantitative results for $S=\nabla$ using the residual based indicators: Endpoint errors (EE) and angular errors (AE) as in \cite{BaScLeRoBlSz:11}, each given by mean and standard deviation, and computational times in seconds.
  \label{tab:OptFlow:nabla}}

\begin{tabular}{ll|ccccc}
\toprule
Benchmark &    Algorithm &     EE-mean &     EE-stddev &     AE-mean &     AE-stddev & Time (s)\\ \hline \hline
Dimetrodon &
    \cite[Algorithm 3]{ours2022semismooth1} without warping & 1.80 & 0.75 & 0.84 & 0.22 & 211.00\\
  & \cite[Algorithm 3]{ours2022semismooth1} with warping  & 0.25 & 0.23 & 0.08 & 0.07 & 2,147.03 \\
  & \cref{alg:optflow_adapt} (with \Code{l2_lagrange}) & 0.41 & 0.37 & 0.13 & 0.12 & 83.96 \\
  &\cref{alg:optflow_adapt} (with \Code{l2-pixel})& 0.42 & 0.38 & 0.14 & 0.12 & 74.62\\
\hline
Grove2 & 
    \cite[Algorithm 3]{ours2022semismooth1} without warping & 2.93 & 0.46 & 1.10 & 0.14 & 285.74 \\
& \cite[Algorithm 3]{ours2022semismooth1} with warping & 1.13 & 1.30 & 0.42 & 0.54 & 1,727.35 \\
& \cref{alg:optflow_adapt} (with l2 lagrange) & 0.40 & 0.56 & 0.11 & 0.16 & 152.76 \\
& \cref{alg:optflow_adapt} (with l2-pixel) & 0.39 & 0.56 & 0.11 & 0.16 & 124.45 \\
  \hline
Grove3 &
    \cite[Algorithm 3]{ours2022semismooth1} without warping & 3.69 & 2.33 & 1.01 & 0.32 & 312.74 \\
 & \cite[Algorithm 3]{ours2022semismooth1} with warping & 1.93 & 2.11 & 0.26 & 0.29 & 1,628.13 \\
& Algorithm 2 (with l2 lagrange) & 1.12 & 1.50 & 0.16 & 0.29 & 168.02 \\
& \cref{alg:optflow_adapt} (with l2-pixel) & 1.10 & 1.49 & 0.16 & 0.29 & 176.56 \\
  \hline
Hydrangea &
    \cite[Algorithm 3]{ours2022semismooth1} without warping & 3.33 & 1.22 & 0.91 & 0.23 & 222.65 \\
&\cite[Algorithm 3]{ours2022semismooth1} with warping & 0.22 & 0.52 & 0.04 & 0.10 & 1,396.16 \\
&\cref{alg:optflow_adapt} (with l2 lagrange) & 0.58 & 0.47 & 0.08 & 0.11 & 106.36 \\
&\cref{alg:optflow_adapt} (with l2-pixel) & 0.56 & 0.44 & 0.07 & 0.11 & 83.33 \\
  \hline
RubberWhale &
    \cite[Algorithm 3]{ours2022semismooth1} without warping & 0.67 & 0.58 & 0.36 & 0.25 & 263.60 \\
&\cite[Algorithm 3]{ours2022semismooth1} with warping & 0.26 & 0.52 & 0.14 & 0.26 & 633.41 \\
&\cref{alg:optflow_adapt} (with l2 lagrange) & 0.37 & 0.59 & 0.20 & 0.33 & 80.61 \\
&\cref{alg:optflow_adapt} (with l2-pixel) & 0.37 & 0.59 & 0.20 & 0.33 & 71.32 \\
  \hline
Urban2 &
    \cite[Algorithm 3]{ours2022semismooth1} without warping & 8.16 & 8.12 & 0.99 & 0.40 & 238.07 \\
&\cite[Algorithm 3]{ours2022semismooth1} with warping & 7.62 & 8.15 & 0.64 & 0.39 & 691.12 \\
&\cref{alg:optflow_adapt} (with l2 lagrange) & 4.53 & 5.32 & 0.27 & 0.30 & 105.81 \\
&\cref{alg:optflow_adapt} (with l2-pixel) & 5.01 & 5.79 & 0.30 & 0.32 & 88.97 \\
  \hline
Urban3 &
    \cite[Algorithm 3]{ours2022semismooth1} without warping & 7.04 & 4.45 & 1.13 & 0.24 & 237.47 \\
&\cite[Algorithm 3]{ours2022semismooth1} with warping & 5.56 & 4.45 & 0.50 & 0.40 & 1,173.47 \\
&\cref{alg:optflow_adapt} (with l2 lagrange) & 1.72 & 2.35 & 0.25 & 0.50 & 137.35 \\
&\cref{alg:optflow_adapt} (with l2-pixel) & 1.78 & 2.35 & 0.25 & 0.50 & 119.87 \\
  \hline
Venus 
 & \cite[Algorithm 3]{ours2022semismooth1} without warping & 3.55 & 1.75 & 1.00 & 0.24 & 125.81 \\
&\cite[Algorithm 3]{ours2022semismooth1} with warping & 1.25 & 1.53 & 0.23 & 0.35 & 1,056.74 \\
&\cref{alg:optflow_adapt} (with l2 lagrange) & 0.72 & 1.13 & 0.15 & 0.34 & 81.04 \\
&\cref{alg:optflow_adapt} (with l2-pixel) & 0.72 & 1.13 & 0.15 & 0.34 & 85.11 \\
  \bottomrule
\end{tabular}
\end{sidewaystable}

For all benchmarks the same manually chosen model parameters were used.
We choose 
$\alpha_1=10$,   
$\alpha_2 = 0$,
$\lambda=1$
to obtain visually pleasing results, cf.\ superiority of L1-TV in
\cite{Dirks:15},
$\beta = 1\cdot 10^{-5}$,
$\gamma_1 = 1\cdot 10^{-4}$,
$\gamma_2 = 1\cdot 10^{-4}$
to balance between speed and quality of the reconstruction and $\vecv{u}_0 := 0$.
For the interior method 
$\epsilon_{\text{newton}} = 1\cdot 10^{-3}$
was chosen.
The mesh is initialized to
$\tfrac{1}{8}$ of the image
resolution, rounded down to integer values, 
$\epsilon_{\text{warp}} = 5\cdot 10^{-2}$
and a constant number of 6
total refinements are carried out before stopping the algorithm.

In \cref{fig:middlebury} we evaluate \cref{alg:optflow_adapt}
visually against the middlebury optical flow benchmark \cite{BaScLeRoBlSz:11} for $S=\nabla$, \nnew{using the residual based indicators in the marking strategy} and using \Code{l2_lagrange} for projecting onto the finite element space.
The color-coded images representing optical flow fields are normalized by the
maximum motion of the ground truth flow data and black areas of the ground
truth data represent unknown flow information, e.g.\ due to occlusion.
A good resemblance of the computed optical flow to the ground truth
and the effect of total variation regularization, i.e.\ sharp edges separating
homogeneous regions, can be seen clearly.
Large displacements are resolved quite well, e.g.\ the fast moving triangle-shaped
object at the bottom of the RubberWhale benchmark, thanks to the warping algorithm
employed.
Using the same example, smaller slow-moving parts adjacent to the larger moving
objects tend to get somewhat distorted however. For quantitatively evaluating our experiments we use the endpoint error and angular error defined as in \cite{BaScLeRoBlSz:11} and summarize our findings in \cref{tab:OptFlow:nabla}. Additionally we also compare in \cref{tab:OptFlow:nabla} the influence of the projections \Code{l2_lagrange} and \Code{l2_pixel} on \cref{alg:optflow_adapt}. In particular, it turns out that independently whether \Code{l2_lagrange} or \Code{l2_pixel} is used, the obtained result are nearly the same and visually hardly distinguishable, see \cref{fig:l2_lagrange:l2_pixel}. 

Remarkable, but maybe not surprisingly, is the speed-up of our adaptive warping technique, which can be clearly seen in \cref{tab:OptFlow:nabla}. Hence \cref{alg:optflow_adapt} not only improves qualitatively the computation of the optical flow of \cite[Algorithm 3]{ours2022semismooth1} but also dramatically shortens the computational time. In particular we see from \cref{tab:OptFlow:nabla} that \cref{alg:optflow_adapt} is even faster than \cite[Algorithm 3]{ours2022semismooth1} without warping and additionally generates a qualitatively much better optical flow. 

\begin{figure}
  \def\imgfraction{0.1428}
  \centering
  \foreach \data in {Dimetrodon,Grove2,Grove3,Hydrangea,RubberWhale,Urban2,Urban3,Venus} {    
    \begin{subfigure}[t]{\imgfraction\columnwidth}
      \includegraphics[width=\textwidth]{fofm_w_c_adaptiveadaptive-warping\data output.png}
    \end{subfigure}%
    \begin{subfigure}[t]{\imgfraction\columnwidth}
      \includegraphics[width=\textwidth]{fofm_w_c_adaptive_pa-w\data output.png}
    \end{subfigure}%
    \begin{subfigure}[t]{\imgfraction\columnwidth}
      \includegraphics[width=\textwidth]{fofm\data ground_truth.png}
    \end{subfigure}%
  }
  \caption{%
    Middlebury Optical Flow Benchmark (Dimetrodon, Grove2, Grove3, Hydrangea,
    RubberWhale, Urban2, Urban3, Venus): \cref{alg:optflow_adapt} with \Code{l2_lagrange} ($1^{\text{st}}$ and $4^{\text{th}}$ column), \cref{alg:optflow_adapt} with \Code{l2_pixel} ($2^{\text{nd}}$ and $5^{\text{th}}$ column) and respective ground truth ($3^{\text{rd}}$ and $6^{\text{th}}$ column)
  }
  \label{fig:l2_lagrange:l2_pixel}
\end{figure}

\begin{figure}
  \centering
  \includegraphics[width=0.49\textwidth]{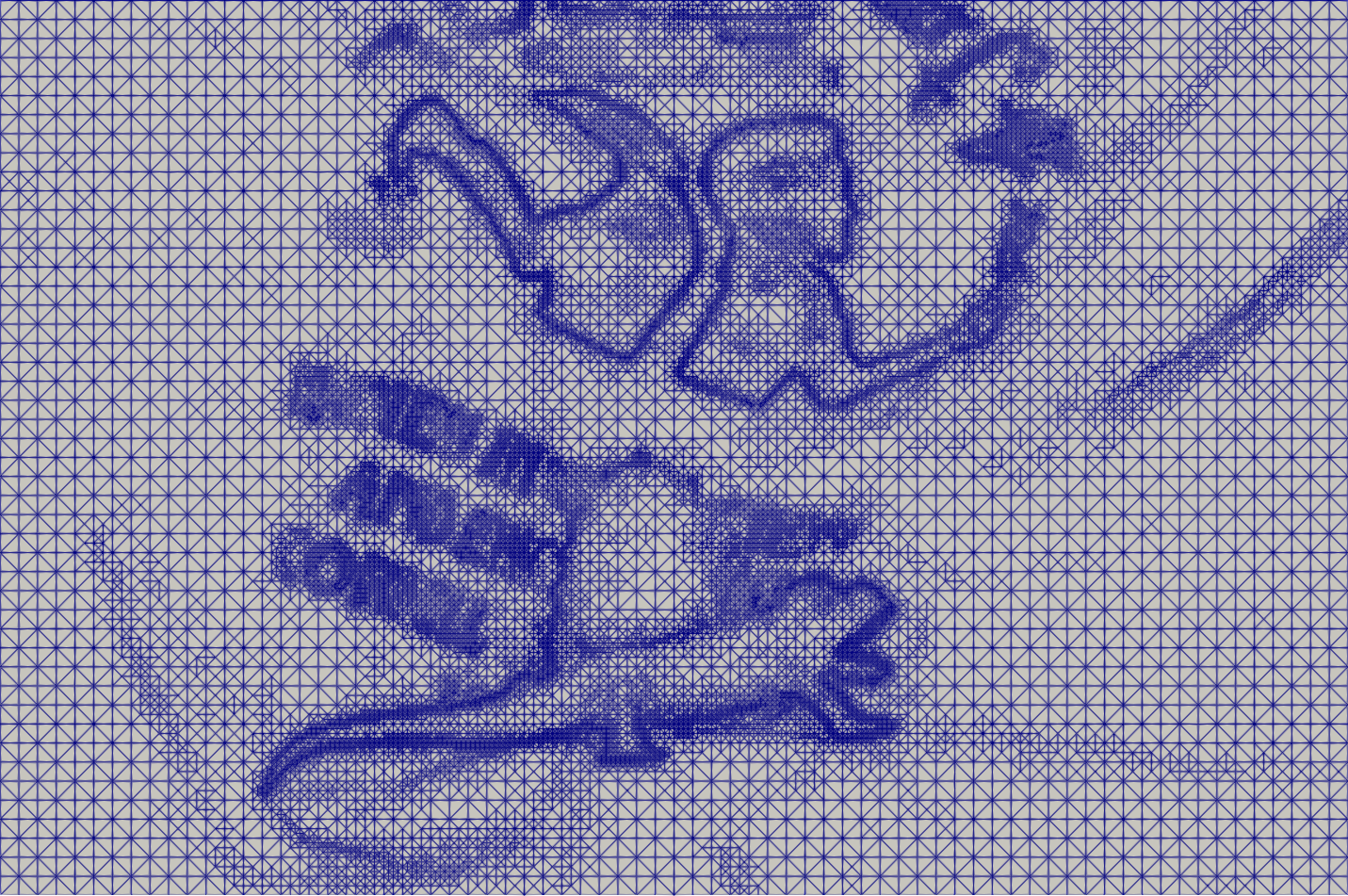}
  \includegraphics[width=0.49\textwidth]{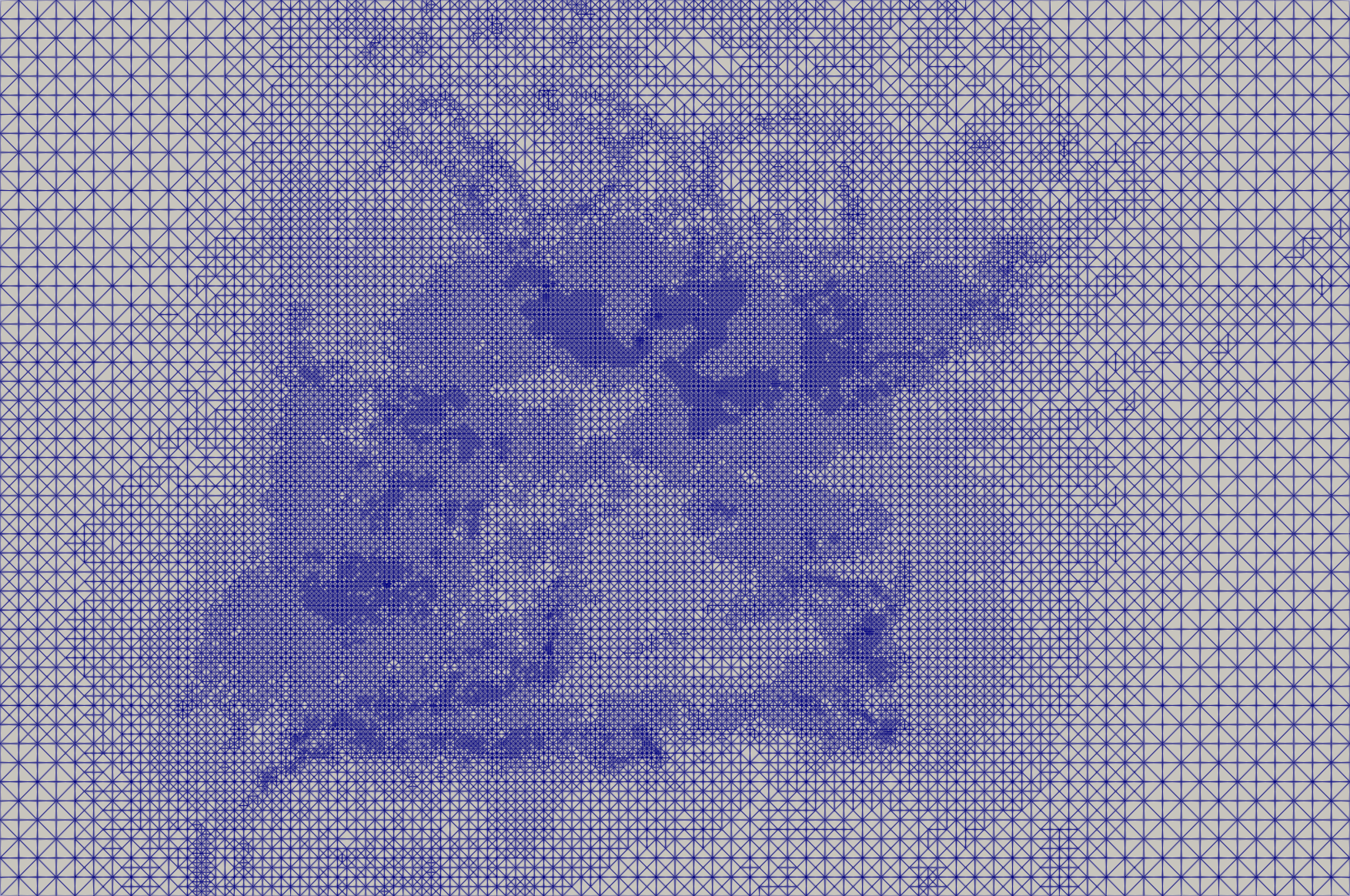}
  \caption{%
    Exemplary adapted mesh for \cref{alg:optflow_adapt} on the Dimetrodon
    Middlebury optical flow benchmark, \nnew{using the residual based indicators (left) and the primal-dual based indicators (right)}. 
  }
  \label{fig:middlebury-mesh}
\end{figure}

It is unclear how much visual improvement more careful or adaptive parameter
selection may give and further study remains to be done.
Exemplary, the adapted mesh for the Dimetrodon benchmark can be seen
in \cref{fig:middlebury-mesh}, where
refinement seems to take place largely around image edges.

\begin{sidewaystable}
\scriptsize
  \centering
\caption{%
    Results for $S=I$ using the residual based indicators: Endpoint errors (EE) and angular errors (AE) as in \cite{BaScLeRoBlSz:11}, each given by mean and standard deviation, and computational time in seconds. \label{tab:OpticalFlow:I}
  }

\begin{tabular}{ll|ccccc}
\toprule
Benchmark &    Algorithm &     EE-mean &     EE-stddev &     AE-mean &     AE-stddev & Time (s) \\ \hline \hline
Dimetrodon &
    \cite[Algorithm 3]{ours2022semismooth1} without warping & 1.81 & 0.75 & 0.84 & 0.22 & 151.10 \\
  & \cite[Algorithm 3]{ours2022semismooth1} with warping  & 0.25 & 0.23 & 0.08 & 0.07 & 1,834.52 \\
  & \cref{alg:optflow_adapt} (with \Code{l2_lagrange}) & 0.43 & 0.40 & 0.14 & 0.13 & 79.52 \\
  & \cref{alg:optflow_adapt} (with \Code{l2-pixel}) & 0.43 & 0.39 & 0.14 & 0.12 & 60.89 \\
\hline
Grove2 & 
    \cite[Algorithm 3]{ours2022semismooth1} without warping & 2.93 & 0.46 & 1.10 & 0.14 & 195.36 \\
  & \cite[Algorithm 3]{ours2022semismooth1} with warping & 1.13 & 1.30 & 0.42 & 0.54 & 1,689.12 \\
  & \cref{alg:optflow_adapt} (with l2 lagrange) & 0.40 & 0.57 & 0.11 & 0.16 & 134.80 \\
  & \cref{alg:optflow_adapt} (with l2-pixel) & 0.39 & 0.56 & 0.11 & 0.16 & 109.99 \\
\hline
Grove3 &
    \cite[Algorithm 3]{ours2022semismooth1} without warping & 3.69 & 2.33 & 1.01 & 0.32 & 211.27 \\
  & \cite[Algorithm 3]{ours2022semismooth1} with warping & 1.93 & 2.10 & 0.27 & 0.29 & 1,579.88 \\
  & Algorithm 2 (with l2 lagrange) & 1.14 & 1.51 & 0.17 & 0.29 & 152.71 \\
  & \cref{alg:optflow_adapt} (with l2-pixel) & 1.12 & 1.49 & 0.16 & 0.29 & 151.63 \\
\hline
Hydrangea &
    \cite[Algorithm 3]{ours2022semismooth1} without warping & 3.33 & 1.22 & 0.91 & 0.23 & 151.16 \\
  & \cite[Algorithm 3]{ours2022semismooth1} with warping & 0.22 & 0.52 & 0.04 & 0.10 & 1,303.31 \\
  & \cref{alg:optflow_adapt} (with l2 lagrange) & 0.60 & 0.50 & 0.08 & 0.12 & 86.75 \\
  & \cref{alg:optflow_adapt} (with l2-pixel) & 0.57 & 0.47 & 0.08 & 0.11 & 72.85 \\
\hline
RubberWhale &
    \cite[Algorithm 3]{ours2022semismooth1} without warping & 0.67 & 0.58 & 0.36 & 0.25 & 187.59 \\
  & \cite[Algorithm 3]{ours2022semismooth1} with warping & 0.26 & 0.52 & 0.14 & 0.26 & 624.25 \\
  & \cref{alg:optflow_adapt} (with l2 lagrange) & 0.39 & 0.59 & 0.21 & 0.33 & 65.65 \\
  & \cref{alg:optflow_adapt} (with l2-pixel) & 0.39 & 0.60 & 0.21 & 0.33 & 57.27 \\
\hline
Urban2 &
    \cite[Algorithm 3]{ours2022semismooth1} without warping & 8.16 & 8.12 & 0.99 & 0.40 & 212.86 \\
  & \cite[Algorithm 3]{ours2022semismooth1} with warping & 7.62 & 8.15 & 0.64 & 0.39 & 637.57 \\
  & \cref{alg:optflow_adapt} (with l2 lagrange) & 4.52 & 5.26 & 0.28 & 0.32 & 83.84 \\
  & \cref{alg:optflow_adapt} (with l2-pixel) & 5.03 & 5.76 & 0.32 & 0.34 & 69.16 \\
\hline
Urban3 &
    \cite[Algorithm 3]{ours2022semismooth1} without warping & 7.04 & 4.45 & 1.13 & 0.24 & 242.47 \\
  & \cite[Algorithm 3]{ours2022semismooth1} with warping & 5.57 & 4.45 & 0.50 & 0.40 & 1,062.80 \\
  & \cref{alg:optflow_adapt} (with l2 lagrange) & 1.76 & 2.34 & 0.25 & 0.50 & 127.76 \\
  & \cref{alg:optflow_adapt} (with l2-pixel) & 1.82 & 2.35 & 0.26 & 0.50 & 100.39 \\
\hline
Venus &
    \cite[Algorithm 3]{ours2022semismooth1} without warping & 3.55 & 1.75 & 1.00 & 0.24 & 121.23 \\
  & \cite[Algorithm 3]{ours2022semismooth1} with warping & 1.25 & 1.52 & 0.23 & 0.35 & 887.40 \\
  & \cref{alg:optflow_adapt} (with l2 lagrange) & 0.73 & 1.13 & 0.15 & 0.34 & 79.46 \\
  & \cref{alg:optflow_adapt} (with l2-pixel) & 0.75 & 1.14 & 0.15 & 0.34 & 80.27 \\
    \bottomrule
\end{tabular}
\end{sidewaystable}

We note that for $S=I$ very similar qualitative and quantitative results are obtained, see \cref{tab:OpticalFlow:I}. A reason for this might be that we choose $\beta = 1 \cdot 10^{-5}$,
which is rather small and hence the choice of $S$ does not play a crucial role for numerically computing a solution. Indeed, as the associated term $\frac{\beta}{2} \|S\vecv{u}\|_{V_S}^2$ is only added for theoretical reasons, one typically wants to choose $\beta$ very small \cite{ours2022semismooth1}.

\nnew{Further, very similar qualitatively results are obtained when the primal-dual based indicators is used in the marking strategy instead of the residual based indicators, see \cref{tab:OpticalFlow-pd}. However, since the indicator needs to be computed on a globally refined mesh in each iteration, the computation times are larger than the ones when the residual based indicators is used. However, it is still in nearly all the experiments faster than \cite[Algorithm 3]{ours2022semismooth1} without warping.
Interestingly, \cref{alg:AFEM} using the primal-dual based indicators generates a different mesh than when the residual based indicators are used, see \cref{fig:middlebury-mesh}. However, the generated optical flows are visually hardly distinguishable, as we observe in the case of Dimetrodon, see \cref{fig:optflow:indicator}.
}
 
\begin{table}
\scriptsize
  \centering
\caption{%
   \nnew{Results for $S=\nabla$ obtained by \cref{alg:optflow_adapt} (with \Code{l2_lagrange}) using the primal-dual based indicators: Endpoint errors (EE) and angular errors (AE) as in \cite{BaScLeRoBlSz:11}, each given by mean and standard deviation, and computational time in seconds.} \label{tab:OpticalFlow-pd}
  }
\begin{tabular}{l|ccccc}
\toprule
{Benchmark} & {EE-mean} & {EE-stddev} & {AE-mean} & {AE-stddev} & {Time (s)} \\
\hline \hline
Dimetrodon   & 0.42 & 0.38 & 0.14 & 0.12 & 139.94 \\
Grove2       & 0.39 & 0.55 & 0.11 & 0.16 & 249.05 \\
Grove3       & 1.15 & 1.50 & 0.17 & 0.29 & 237.76 \\
Hydrangea    & 0.52 & 0.44 & 0.07 & 0.10 & 151.24 \\
RubberWhale  & 0.39 & 0.60 & 0.21 & 0.33 & 140.42 \\
Urban2       & 4.38 & 5.16 & 0.29 & 0.38 & 175.19 \\
Urban3       & 1.76 & 2.35 & 0.25 & 0.50 & 248.80 \\
Venus        & 0.70 & 1.11 & 0.14 & 0.34 & 154.45 \\
\bottomrule
\end{tabular}
  
\end{table}

\begin{figure}
  \def\imgfraction{0.49}
  \centering
  \foreach \data in {Dimetrodon} {    
    \begin{subfigure}[t]{\imgfraction\columnwidth}
      \includegraphics[width=\textwidth]{fofm_w_c_adaptiveadaptive-warping\data output.png}
    \end{subfigure}\hfill%
    \begin{subfigure}[t]{\imgfraction\columnwidth}
      \includegraphics[width=\textwidth]{fofm_w_c_adaptiveadaptive-warping-pd\data output.png}
    \end{subfigure}
  }
  \caption{%
  \cref{alg:optflow_adapt} with \Code{l2_lagrange} applied on Dimetrodon: residual based indicators (left), primal-dual based indicators (right)
  }
  \label{fig:optflow:indicator}
\end{figure}

\backmatter

%
%
%

\section*{Declarations}

\bmhead{Funding}
This work was partly supported by the Ministerium f\"ur Wissenschaft, Forschung und Kunst of Baden-W\"urttemberg (Az: 7533.-30-10/56/1) through the RISC-project ``Automatische Erkennung von bewegten Objekten in hochaufl\"osenden Bildsequenzen mittels neuer Gebietszerlegungsverfahren'', by the Deutsche Forschungsgemeinschaft (DFG, German Research Foundation) under Germany's Excellence Strategy -- EXC 2075 -- 390740016, and by the Crafoord Foundation through the project ``Locally Adaptive Methods for Free Discontinuity Problems''.

\bmhead{Conflict of interest}
The authors have no competing interests to declare that are relevant to the content of this article.

\bmhead{Code availability}
The code is available as a Julia package and can be found at \url{https://gitlab.mathematik.uni-stuttgart.de/stephan.hilb/SemiSmoothNewton.jl}.

%
%

\begin{appendices}
\renewcommand{\theequation}{\thesection.\arabic{equation}}
\section{Proof of \cref{Thm:ErrorEstimate}}\label{Sec:ProofErrorEstimate}
Denote $\vecv{e}_h := \vecv{u}_{h} - \vecv{u} \in V$ and bound using
\cref{lem:energy-strong-convexity}:
\begin{align*}
  \tfrac{1}{2} \|\vecv{u}_{h} - \vecv{u}\|_B^2
  &\le E(\vecv{u}_{h}) - E(\vecv{u}) \\
  &= F(\vecv{u}_{h}) - F(\vecv{u}) + G(\Lambda \vecv{u}_{h}) - G(\Lambda \vecv{u})\\
  &\le \langle \partial F(\vecv{u}_{h}), \vecv{e}_h\rangle_{V^*,V} + \langle \partial G(\Lambda \vecv{u}_{h}), \Lambda \vecv{e}_h\rangle_{L^2} \\
  &= a_B(\vecv{u}_{h}, \vecv{e}_h) - l(\vecv{e}_h) + \langle \partial G(\Lambda \vecv{u}_{h}), \Lambda
  \vecv{e}_h\rangle_{L^2}.
\end{align*}
Depending on $S$, we will now 
further bound
the error $\tfrac{1}{2}\|\vecv{u}_{h} - \vecv{u}\|_B^2$.

\subsubsection*{\Cref{setting-nabla}: the Case ${S = \nabla}$}

We will make use of the fact that $a_B$ is coercive on $V$, see
\cref{assumption-a-coercivity},
to derive an estimate on the error
$\|\vecv{e}_h\|_{H^1(\Omega)^m}$.

Due to \eqref{eq:opt_cond_discrete_conformity} and the discrete
optimality condition \eqref{eq:opt_cond_discrete_f} we have for any bounded linear
operator $I_h: V \to V_h$:
\begin{equation}
  \begin{aligned} \label{eq:apost_prelim}
    \tfrac{1}{2} \|\vecv{u}_{h} - \vecv{u}\|_B^2
    &\le a_B(\vecv{u}_{h}, \vecv{e}_h) - l(\vecv{e}_h) + \<\vec{p}_{h}, \Lambda \vecv{e}_h\>_{L^2} \\
    &= a_B(\vecv{u}_{h}, \vecv{e}_h - I_h \vecv{e}_h) - l(\vecv{e}_h - I_h \vecv{e}_h) + \<\vec{p}_{h}, \Lambda (\vecv{e}_h - I_h
    \vecv{e}_h)\>_{L^2} \\
    &= \<\alpha_2 (T\vecv{u}_{h} - g) + p_{h,1}, T(\vecv{e}_h - I_h \vecv{e}_h)\>_{L^2} \\
      &\qquad + \<\beta S\vecv{u}_{h}, S(\vecv{e}_h - I_h \vecv{e}_h)\>_{L^2}
      + \<\vec{p}_{h,2}, \nabla (\vecv{e}_h - I_h \vecv{e}_h)\>_{L^2}.
  \end{aligned}
\end{equation}
Actually we choose $I_h: V \to V_h$ to be a quasi-interpolation operator which satisfies
the interpolation estimates \cite[Proposition 1.3]{Verfurth:13} \cite[Theorem 1.7]{AinsworthOden:00}:
\begin{align}
  \|\vecv{v} - I_h \vecv{v}\|_{L^2(K)}
    &\leq c_1 h_K \|\nabla \vecv{v}\|_{L^2(\omega_K)}, \label{eq:interpolation-estimate-1} \\
  \|\vecv{v} - I_h \vecv{v} \|_{L^2(F)}
    &\leq c_2 h_F^{\frac{1}{2}}\|\nabla \vecv{v}\|_{L^2(\omega_F)}, \label{eq:interpolation-estimate-2}
\end{align}
where $c_1,c_2>0$, $\omega_K$, $\omega_F$ denote the union of triangles which share a common
vertex with the cell $K$ or the facet $F$ respectively and $h_K$, $h_F$ denote
the diameter of $K$ and $F$ respectively.

Utilizing Green's formula, see e.g.\ \cite[Proposition 2.4]{Necas:2011}, in \eqref{eq:apost_prelim} yields 
\begin{align*}
  \tfrac{1}{2} &\|\vecv{u}_{h} - \vecv{u}\|_B^2
  \le \big\<\alpha_2(T\vecv{u}_{h} - g) + p_{h,1},
      T(\vecv{e}_h - I_h \vecv{e}_h)\big\>_{L^2} \\
      &\phantom{\|\vecv{u}_{h} - \vecv{u}\|_B^2
  \le }+
    \big\<\beta \nabla \vecv{u}_{h} + \vec{p}_{h,2},
      \nabla (\vecv{e}_h - I_h \vecv{e}_h) \big\>_{L^2} \\
  &= \sum_{K\in \mathcal T}\bigg( \Big\<\alpha_2 T^*(T\vecv{u}_{h} -g) + T^*p_{h,1} - \beta
    \Laplace \vecv{u}_{h} - \Div \vec{p}_{h,2}, \vecv{e}_h - I_h\vecv{e}_h\Big\>_{L^2(K)^m} \\
  &\qquad\qquad + \Big\<\vecv{n}^\top(\beta \nabla \vecv{u}_{h} + \vec{p}_{h,2}), \vecv{e}_h -
    I_h\vecv{e}_h\Big\>_{L^2(\partial K)^m} \bigg).
\end{align*}
Denoting $\vecv{\xi} := \alpha_2 T^*(T\vecv{u}_{h} -g) + T^* p_{h,1} - \beta \Laplace \vecv{u}_{h} -
\Div \vec{p}_{h,2} \in L^2(\Omega)^m$ and $\vec{\zeta} := \beta \nabla \vecv{u}_{h} + \vec{p}_{h,2} \in
L^2(\Omega)^{d\times m}$ we obtain using interpolation estimates
\eqref{eq:interpolation-estimate-1,eq:interpolation-estimate-2}:
\begin{align*}
  \tfrac{1}{2} \|\vecv{u}_{h} - \vecv{u}\|_B^2
  &\le \sum_{K \in \mathcal T} \Big(\<\vecv{\xi}, \vecv{e}_h - I\vecv{e}_h\>_{L^2(K)^m} + \<\vec{n}^\transpose \vec{\zeta}, \vecv{e}_h - I\vecv{e}_h\>_{L^2(\partial K)^m} \Big) \\
  &= \sum_{K \in \mathcal T} \<\vecv{\xi}, \vecv{e}_h - I\vecv{e}_h\>_{L^2(K)^m} + \sum_{F\in\Gamma} \<[\vec{n}^\transpose\vec{\zeta}], \vecv{e}_h - I\vecv{e}_h\>_{L^2(F)^m} \\
  &\le \max\{c_1,c_2\} \Big( \sum_{K\in\mathcal T} \|\nabla \vecv{e}_h\|_{L^2(\omega_K)^{d\times m}} h_K\|\vecv{\xi}\|_{L^2(K)^m} \\
  &\qquad\qquad\qquad\quad+ \sum_{F\in\Gamma} \|\nabla \vecv{e}_h\|_{L^2(\omega_F)^{d\times m}} h_F^{\frac{1}{2}} \|[\vec{n}^\transpose\vec{\zeta}]\|_{L^2(F)^m} \Big) \\
  &\le \max\{c_1,c_2\}\Big( \sum_{K\in\mathcal T} \|\nabla \vecv{e}_h\|_{L^2(\omega_K)^{d\times m}}^2  + \sum_{F\in\Gamma} \|\nabla \vecv{e}_h\|_{L^2(\omega_F)^{d\times m}}^2\Big)^{\frac{1}{2}} \\
  &\qquad\qquad\qquad\quad \Big(
    \sum_{K\in\mathcal T} h_K^2\|\vecv{\xi}\|_{L^2(K)^m}^2 +
    \sum_{F\in\Gamma}h_F \|[\vec{n}^\transpose\vec{\zeta}]\|_{L^2(F)^m}^2 \Big)^{\frac{1}{2}}.
\end{align*}
Let $B(K)$ be the set of facets of $K$ belonging to $\partial\Omega$ and $I(K)$ all other facets of $K$. 
Realizing that $\omega_F \subset \omega_K$ for all to $F$ adjacent cells $K$ and having at most one boundary facet at each cell
we may bound 
\begin{equation}\label{Eq:cell_inequ}
\begin{split}
  \sum_{F\in\Gamma} \|&\nabla \vecv{e}_h\|_{L^2(\omega_F)^{d\times m}}^2 \\
  &=
  \sum_{K\in\mathcal T} \left(\sum_{F\in B(K)}\|\nabla \vecv{e}_h\|_{L^2(\omega_F)^{d\times m}}^2 + \frac12 \sum_{F \in I(K)}\|\nabla \vecv{e}_h\|_{L^2(\omega_F)^{d\times m}}^2\right) \\
  &\le
  \sum_{K\in\mathcal T} \left(\sum_{F\in B(K)}\|\nabla \vecv{e}_h\|_{L^2(\omega_K)^{d\times m}}^2 + \frac12 \sum_{F \in I(K)}\|\nabla \vecv{e}_h\|_{L^2(\omega_K)^{d\times m}}^2\right) \\
  &\le 2 \sum_{K\in\mathcal T} \|\nabla \vecv{e}_h\|_{L^2(\omega_K)^{d\times m}}^2
  \le 2 c_{\mathcal T}^2 \|\nabla \vecv{e}_h\|_{L^2(\Omega)^{d\times m}}^2,
\end{split}
\end{equation}
where $c_{\mathcal T}$ denotes the shape regularity constant of the mesh independent of the mesh size.
We thus arrive at
\begin{align*}
  &\tfrac{1}{2} \|\vecv{u}_{h} - \vecv{u}\|_B^2 \\
  &\quad\le \max\{c_1,c_2\}\Big( 3 \sum_{K\in\mathcal T} \|\nabla \vecv{e}_h\|_{L^2(\omega_K)^{d\times m}}^2 \Big)^{\frac{1}{2}}\Big(\sum_{K\in\mathcal T} h_K^2\|\vecv{\xi}\|_{L^2(K)^m}^2\\
  &\phantom{\hspace{7.5cm}}  + \sum_{F\in\Gamma}h_F \|[\vec{n}^\transpose\vec{\zeta}]\|_{L^2(F)^m}^2 \Big)^{\frac{1}{2}} \\
  &\quad\le \sqrt{3} \max\{c_1,c_2\} c_{\mathcal T} \|\nabla \vecv{e}_h\|_{L^2(\Omega)^{d\times m}}\Big(\sum_{K\in\mathcal T} h_K^2\|\vecv{\xi}\|_{L^2(K)^m}^2\\
  &\phantom{\hspace{7.5cm}} + \sum_{F\in\Gamma}h_F \|[\vec{n}^\transpose\vec{\zeta}]\|_{L^2(F)^m}^2 \Big)^{\frac{1}{2}} \\
  &\quad\le \sqrt{3} \max\{c_1,c_2\} c_{\mathcal T}
  \|\vecv{e}_h\|_{H^1(\Omega)^m}\Big(\sum_{K\in\mathcal T} h_K^2\|\vecv{\xi}\|_{L^2(K)^m}^2 \\
  &\phantom{\hspace{7.5cm}}+ \sum_{F\in\Gamma}h_F \|[\vec{n}^\transpose\vec{\zeta}]\|_{L^2(F)^m}^2 \Big)^{\frac{1}{2}}.
\end{align*}
Together with coercivity $\|\vecv{v}\|_B^2 = a_B(\vecv{v},\vecv{v}) \ge c_B \|\vecv{v}\|_V^2 =
c_B \|\vecv{v}\|_{H^1(\Omega)^m}^2$ from \cref{assumption-a-coercivity} we arrive at an a-posteriori error bound
\begin{align*}
  \|\vecv{u}_{h} - \vecv{u}\|_{H^1(\Omega)^m}^2 &\le C \bigg(\sum_{K \in \mathcal T} h_K^2 \|\alpha_2 T^*(T\vecv{u}_{h} - g) +
    T^*p_{h,1}
    - \beta \Laplace \vecv{u}_{h} - \Div \vec{p}_{h,2}\|_{L^2(K)^m}^2 \\
  &\qquad\qquad + \sum_{F\in\Gamma} h_F \|[\vec{n}^\transpose(\beta\nabla \vecv{u}_{h} + \vec{p}_{h,2})]\|_{L^2(F)^m}^2 \bigg)
\end{align*}
with constant $C := \frac{12}{c_B} \max\{c_1^2,c_2^2\}c_{\mathcal T}^2$.

\subsubsection*{\Cref{setting-identity}: the case ${S = I}$}

Using the fact that $a_B$ is coercive in this setting with respect to $L^2$, see \cref{assumption-a-coercivity}, we derive an estimate for the error
$\|\vecv{e}_h\|_{L^2(\Omega)^m}$.

Let $I_h: L^2(\Omega)^m \to V_h$  be a bounded and linear quasi-interpolation operator satisfying the interpolation estimate \cite[Theorem
1.7]{AinsworthOden:00}:
\begin{align}
  \|\vecv{v} - I_h \vecv{v} \|_{L^2(F)}
    \leq c_3 h_F^{-\frac{1}{2}}\|\vecv{v}\|_{L^2(\omega_F)}, \label{eq:interpolation-estimate-3}
\end{align}
where $c_3 > 0$. Then, as above in \cref{setting-nabla}, cf.\ \eqref{eq:apost_prelim}, we obtain from \eqref{eq:opt_cond_discrete_conformity} and the discrete
optimality condition \eqref{eq:opt_cond_discrete_f} that we have
\begin{equation}
  \begin{aligned} \label{eq:apost_prelim2}
    \tfrac{1}{2} \|\vecv{u}_{h} - \vecv{u}\|_B^2
    &= \<\alpha_2 (T\vecv{u}_{h} - g) + p_{h,1}, T(\vecv{e}_h - I_h \vecv{e}_h)\>_{L^2} \\
      &\qquad + \<\beta S\vecv{u}_{h}, S(\vecv{e}_h - I_h \vecv{e}_h)\>_{L^2}
      + \<\vec{p}_{h,2}, \nabla (\vecv{e}_h - I_h \vecv{e}_h)\>_{L^2}.
  \end{aligned}
\end{equation}
Applying Green's formula to \eqref{eq:apost_prelim2} yields the estimate
\begin{align*}
  \tfrac{1}{2} a_B(&\vecv{e}_h, \vecv{e}_h)
  \le \<\alpha_2 T^*(T \vecv{u}_{h} - g)
      + T^*p_{h,1} + \beta \vecv{u}_{h}, \vecv{e}_h - I_h\vecv{e}_h\>_{L^2}
      \\
      &\phantom{\vecv{e}_h, \vecv{e}_h)
  \le \<}+ \<\vec{p}_{h,2}, \nabla (\vecv{e}_h - I_h\vecv{e}_h) \>_{L^2} \\
  &\le \sum_{K \in \mathcal T} \bigg( \Big\<
      \alpha_2 T^*(T \vecv{u}_{h} -g) + T^* p_{h,1}
      + \beta \vecv{u}_{h} - \Div \vec{p}_{h,2}, \vecv{e}_h - I_h \vecv{e}_h\Big\>_{L^2(K)^m} \\
    &\phantom{\vecv{e}_h, \vecv{e}_h)
  \le \<} + \Big\<\vec{n}^\transpose\vec{p}_{h,2}, \vecv{e}_h - I_h \vecv{e}_h \Big\>_{L^2(\partial K)^m}
    \bigg).
\end{align*}
Denoting $\vecv{\xi} := \alpha_2 T^*(T \vecv{u}_{h} -g) + T^* p_{h,1} + \beta \vecv{u}_{h} - \Div
\vec{p}_{h,2}$ and $\vec{\zeta}:= \vec{n}^\transpose \vec{p}_{h,2}$ for brevity one continues to derive using
\eqref{eq:interpolation-estimate-3}, Cauchy-Schwarz inequality and \eqref{Eq:cell_inequ}: 
\begin{align*}
  \tfrac{1}{2} &a_B(\vecv{e}_h, \vecv{e}_h) \\
  &\le \sum_{K \in \mathcal T} \|\vecv{\xi}\|_{L^2(K)^m} \|\vecv{e}_h - I_h \vecv{e}_h\|_{L^2(K)^m}
    + \sum_{F\in\Gamma} \|[\vec{\zeta}]\|_{L^2(F)^m} \|\vecv{e}_h - I_h \vecv{e}_h\|_{L^2(F)^m} \\
  &\stackrel[]{\eqref{eq:interpolation-estimate-3}}{\le} \sum_{K \in \mathcal T} \|\vecv{\xi}\|_{L^2(K)^m} \|\vecv{e}_h - I_h \vecv{e}_h\|_{L^2(K)^m}
    + \sum_{F\in\Gamma} \|[\vec{\zeta}]\|_{L^2(F)^m} h_F^{-1/2} c_3 \|\vecv{e}_h\|_{L^2(\omega_F)^m} \\ 
 &\stackrel[]{\substack{\text{C.-S.}\\\text{ineq.}}}{\le} \left(\sum_{K \in \mathcal T} \|\vecv{e}_h - I_h \vecv{e}_h\|_{L^2(K)^m}^2 + \sum_{F\in\Gamma} c_3^2 \|\vecv{e}_h\|_{L^2(\omega_F)^m}^2 \right)^{\frac{1}{2}} \Bigg( \sum_{K \in \mathcal T} \|\vecv{\xi}\|_{L^2(K)^m}^2 \\
 &\phantom{\hspace{8cm}}
       + \sum_{F\in\Gamma} h_F^{-1} \|[\vec{\zeta}]\|_{L^2(F)^m}^2  \Bigg)^{\frac{1}{2}}\\
 & 
 \le
 \left(\sum_{K \in \mathcal T} \tilde{c} \|\vecv{e}_h\|_{L^2(K)^m}^2 + \sum_{F\in\Gamma} c_3^2 \|\vecv{e}_h\|_{L^2(\omega_F)^m}^2 \right)^{\frac{1}{2}} \Bigg( \sum_{K \in \mathcal T} \|\vecv{\xi}\|_{L^2(K)^m}^2\\
 &\phantom{\hspace{8cm}}   + \sum_{F\in\Gamma} h_F^{-1} \|[\vec{\zeta}]\|_{L^2(F)^m}^2  \Bigg)^{\frac{1}{2}}\\
 &\stackrel[]{\eqref{Eq:cell_inequ}}{\le}
 \left(\sum_{K \in \mathcal T} \tilde{c} \|\vecv{e}_h\|_{L^2(\omega_K)^m}^2 + 2 c_3^2 \|\vecv{e}_h\|_{L^2(\omega_K)^m}^2 \right)^{\frac{1}{2}} \Bigg( \sum_{K \in \mathcal T} \|\vecv{\xi}\|_{L^2(K)^m}^2\\
 &\phantom{\hspace{8cm}}  + \sum_{F\in\Gamma} h_F^{-1} \|[\vec{\zeta}]\|_{L^2(F)^m}^2  \Bigg)^{\frac{1}{2}}\\
 &\stackrel[]{\eqref{Eq:cell_inequ}}{\le} \left(\tilde{c} + 2 c_3^2\right)^{\frac{1}{2}}
 c_{\mathcal T} \|\vecv{e}_h\|_{L^2(\Omega)^m} \left( \sum_{K \in \mathcal T} \|\vecv{\xi}\|_{L^2(K)^m}^2 + \sum_{F\in\Gamma} h_F^{-1} \|[\vec{\zeta}]\|_{L^2(F)^m}^2  \right)^{\frac{1}{2}},
\end{align*}
where $\tilde{c} > 0$. 
Together with the $L^2$-coercivity of $a_B$ from \cref{assumption-a-coercivity} we arrive at
\begin{align*}
  \|\vecv{u}_{h} - \vecv{u}\|_{L^2(\Omega)^m}^2
  &\le C \bigg( \sum_{K \in \mathcal T} \|\alpha_2 T^*(T\vecv{u}_{h} - g) + T^* p_{h,1} +
    \beta \vecv{u}_{h} - \Div \vec{p}_{h,2}\|_{L^2(K)^m}^2  \\
  &\qquad\qquad + \sum_{F\in\Gamma} h_F^{-1} \|[\vec{n}^\transpose\vec{p}_{h,2}]\|_{L^2(F)^m}^2 \bigg)
\end{align*}
for some constant $C > 0$.

\end{appendices}


\end{document}